\documentclass[a4paper]{article}
\pdfoutput=1

\usepackage{graphicx,wrapfig}
\usepackage{lipsum}
\usepackage{amsmath,amssymb,mathrsfs}
\usepackage{hyperref}
\usepackage{url}
\usepackage{mathtools}
\usepackage{changepage}
\usepackage{multirow}
\usepackage{wrapfig}
\usepackage{subfig}
\usepackage{color}
\usepackage{epsfig,enumerate,amsmath,amsfonts,amssymb,amsthm,mathrsfs,ifpdf}
\usepackage{indentfirst,relsize}
\usepackage{setspace,graphicx}
\usepackage{latexsym}
\usepackage[all]{xy}
\usepackage[usenames,dvipsnames]{pstricks}
\usepackage{pst-grad} 
\usepackage{pst-plot} 
\usepackage[margin = 2.50cm]{geometry}

\usepackage{algorithm}

\usepackage[noend]{algpseudocode}

\def\Z{\ensuremath{\mathbb{Z}}}

\newcommand{\remove}[1]{}

\newtheorem{theo}{Theorem}
\newtheorem{lem}[theo]{Lemma}
\newtheorem{pre}[theo]{Proposition}

\newtheorem{cl}[theo]{Claim}

\newtheorem{obs}[theo]{Observation}

\newtheorem{definition}{Definition}[section]
\graphicspath{{./images/}}

\title{{Liar's vertex-edge domination in subclasses of chordal graphs}}

\author{Debojyoti Bhattacharya\footnote{Indian Instittute of Technology Patna, Bihta, 801106, Bihar, India. email: debojyoti\_2021ma11@iitp.ac.in} \and Subhabrata Paul\footnote{Indian Instittute of Technology Patna, Bihta, 801106, Bihar, India. email: subhabrata@iitp.ac.in}}

\date{}
\begin{document}
\maketitle
\begin{abstract}
	Let $G=(V, E)$ be an undirected graph. The set $N_G[x]=\{y\in V|xy\in E\}\cup \{x\}$ is called the closed neighbourhood of a vertex $x\in V$ and for an edge $e=xy\in E$, the closed neighbourhood of $e$ is the set $N_G[x]\cup N_G[y]$, which is denoted by $N_G[e]$ or $N_G[xy]$. A set $L\subseteq V$ is called \emph{liar's vertex-edge dominating set} of a graph $G=(V,E)$ if for every $e_i\in E$, $|N_G[e_i]\cap L|\geq 2$ and for every pair of distinct edges $e_i,e_j\in E$, $|(N_G[e_i]\cup N_G[e_j])\cap L|\geq 3$. The notion of liar's vertex-edge domination arises naturally from some applications in communication networks. Given a graph $G$, the \textsc{Minimum Liar's Vertex-Edge Domination Problem} (\textsc{MinLVEDP}) asks to find a liar's vertex-edge dominating set of $G$ of minimum cardinality. In this paper, we study this problem from an algorithmic point of view. We design two linear time algorithms for \textsc{MinLVEDP} in block graphs and proper interval graphs, respectively. On the negative side, we show that the decision version of liar's vertex-edge domination problem is NP-complete for undirected path graphs.
\end{abstract}



\textbf{keyword}
	Liar's vertex-edge dominating set, NP-completeness, block graphs, proper interval graphs, undirected path graphs

\section{Introduction}
Let $G=(V, E)$ be a graph with vertex set $V$ and edge set $E$. The open neighbourhood and closed neighbourhood of a vertex $v\in V$ is defined by the sets $N_G(v)=\{u\in V|uv\in E\}$ and $N_G[v]=N_G(v)\cup \{v\}$, respectively. The closed neighbourhood of an edge $e=xy\in E$ is the set $N_G[e]=N_G[xy]=N_G[x]\cup N_G[y]$. An edge $e$ is said to be \emph{vertex-edge dominated} (or \emph{ve-dominated}) by a vertex $v$ if $v\in N_G[e]$. A subset $D\subseteq V$ is called a \emph{vertex-edge dominating set} (or a \emph{ve-dominating set}) of $G$ if for every edge $e\in E$, $N_G[e]$ contains at least one vertex of $D$. The cardinality of a minimum ve-dominating set is called the \emph{ve-domination number} of $G$. Peters introduced vertex-edge domination in his Ph.D. thesis in 1986 \cite{peters}. But it did not get much attention until Lewis proved some interesting results \cite{lewis}. Since then, vertex-edge domination has been extensively studied in literature \cite{boutrig2016vertex,jena,paul2,zylinski2019vertex}. For some non negative integer $k$, an edge $e$ is \emph{$k$-vertex edge dominated}(or \emph{$k$-ve dominated}) by $D_k\subseteq V$ if $D_k$ contains at least $k$ vertices of $N_G[e]$. A subset $D_k\subseteq V$ is called a \emph{$k$-ve dominating set} of $G$ if every edge $e\in E$ is $k$-ve dominated by $D_k$, for some integer $k\geq 1$. The cardinality of a minimum $k$-ve dominating set is called the \emph{$k$-ve domination number} of $G$. This variation is also well studied in literature \cite{krishna,li2023polynomial,naresh}. Different variations of vertex-edge domination, namely independent ve-domination \cite{lewis,paul}, total vertex-edge domination \cite{ahangar2021total,Totalve-domchellali}, global vertex-edge domination \cite{chitra2012global,globalvedom} etc are studied in the literature. 

In this article, we study a variation of the vertex-edge domination problem, namely liar's vertex-edge domination problem, which arises naturally from the following application:

Let us consider a graph $G=(V, E)$ representing a communication network where each vertex $v\in V$ is a communication device and each edge $e=uv$ is a communication channel between devices $u$ and $v$. To ensure the smooth functioning of this network, every channel needs to be constantly monitored so that any broken or damaged channel can be identified immediately. This monitoring can be done by sentinels placed at vertices of $G$. The job of a sentinel placed at a vertex $v$ is to detect the broken or damaged channel incident to $N_G[v]$, and then report that broken or damaged channel. To optimize the cost, it is required to find the minimum set of vertices at which the sentinels are to be placed. Note that a minimum ve-dominating set of $G$ is a solution to that optimization problem. If any sentinel fails to detect a broken or damaged channel, then to correctly identify the broken or damaged channel, each channel is required to be monitored by more than one sentinel. A minimum $2$-ve dominating set of $G$ is a solution to this optimization problem. Now, it might happen that a sentinel placed at $v$ correctly detects the broken or damaged channel incident to a vertex of $N_G[v]$ but misreports an undamaged channel incident to a vertex of $N_G[v]$ as a damaged channel. Also, let us assume that every sentinel in the neighbourhood of a broken or damaged edge $e$ detects $e$ as a damaged edge but at most one sentinel at $v$ in $N_G[e]$ is misreporting an undamaged channel incident to a vertex of $N_G[v]$ as a damaged channel. Clearly, a $2$-ve dominating set is not sufficient to correctly detect the broken or damaged edge. To deal with such problems \emph{liar's vertex-edge dominating set} (or liar's ve-dominating set) was introduced and the following definition was established in \cite{liar'svefirst}.

\begin{definition}
	A vertex set $L\subseteq V$ is \emph{liar's ve-dominating set} of a graph $G=(V,E)$ if 
	\begin{itemize}
		\item[(i)] for every $e_i\in E$, $|N_G[e_i]\cap L|\geq 2$ and
		\item[(ii)] for every pair of distinct edges $e_i$ and $e_j$, $|(N_G[e_i]\cup N_G[e_j])\cap L|\geq 3$. 
	\end{itemize}
	The minimum cardinality of a liar's ve-dominating set of $G$ is called the \emph{liar's ve-domination number}. 
\end{definition}

In \cite{liar'svefirst}, authors have proven that the \textsc{DecideLVEDP} is NP-complete for chordal and bipartite graphs but \textsc{MinLVEDP} can be solved in linear time in trees. In addition, they have designed two approximation algorithms for \textsc{MinLVEDP} in general graphs and $p$-claw free graphs and proved that \textsc{MinLVEDP} can not be approximated within $\frac{1}{2}(\frac{1}{8}-\epsilon)\ln |V|$, for any $\epsilon>0$. Also, they have proved that the \textsc{MinLVEDP} is APX-complete in bounded degree graphs and $p$-claw free graphs for $p\geq 6$.

In this paper, we have studied the algorithmic status of liar's ve-domination problem for subclasses of chordal graphs. The minimum liar's ve-domination problem and its corresponding decision version are described as follows:

\noindent\underline{\textsc{Minimum Liar's VE-Domination Problem} (\textsc{MinLVEDP})}

\noindent\emph{Instance}: A graph $G=(V,E)$.

\noindent\emph{Output}: A minimum liar's ve-dominating set $L_{ve}$ of $G$.

\noindent\underline{\textsc{Liar's VE-Domination Decision Problem} (\textsc{DecideLVEDP})}

\noindent\emph{Instance}: A graph $G=(V,E)$ and an integer $k$.

\noindent\emph{Question}: Does there exist a liar's ve-dominating set of size at most $k$?

The results presented in this paper are summarized as follows:
\begin{itemize}
	\item[(i)]  We propose a linear time algorithm for \textsc{MinLVEDP} in block graphs.
	\item[(ii)] We design a linear time algorithm for \textsc{MinLVEDP} in proper interval graphs.
	\item[(ii)] We prove that \textsc{DecideLVEDP} is NP-complete for undirected path graphs.
\end{itemize}

The rest of the paper is organized as follows. In Section $2$, we discuss the necessary preliminaries. After that, in Section $3$, we propose a linear time algorithm for \textsc{MinLVEDP} in block graphs and a linear time algorithm for \textsc{MinLVEDP} in proper interval graphs.
In Section $4$, we prove the NP-completeness result for \textsc{DecideLVEDP} in undirected path graphs. Finally, Section $5$ concludes the paper with some directions for further study.

\section{Preliminaries}

A vertex $v$ is called a cut vertex of a graph $G$ if removing $v$ from $G$ increases the number of components in $G$. For a graph $G$, the maximal connected subgraph of $G$ without any cut vertex is called a block. A graph $G$ is called a block graph if every block of $G$ is a complete subgraph. The intersection of two blocks contains at most one cut vertex. Two blocks are adjacent if their intersection contains a cut vertex of $G$. A block containing only one cut vertex is called an end block of $G$. A cut-tree $T_G$ of a block graph $G$ is the tree whose vertices are the blocks and cut vertices of $G$ and every vertex of the tree corresponding to the blocks of $G$ has an edge with every vertex corresponding to the cut vertices of that block. A block graph $G$ is called \emph{star block} if the corresponding cut-tree $T_G$ is a star. 

A graph $G=(V, E)$ is called an intersection graph of a family of subsets $\mathcal{F}$ if every vertex of $G$ corresponds to a subset in $\mathcal{F}$ and there is an edge between two vertices if the intersection of the corresponding subsets of $\mathcal{F}$ is non-empty. A graph $G=(V, E)$ is called an \emph{interval graph} if it is the intersection graph of some intervals on the real line. In other words, corresponding to the vertices of $G$, there exists a set of intervals $\mathcal{I}=\{I(u): ~ u\in V \}$ such that there is an edge between two vertices $u$ and $v$ in $G$ if and only if their corresponding intervals ($I(u)$ and $I(v)$, respectively) intersect. The set of intervals $\mathcal{I}$ is called the \emph{interval representation} of the interval graph $G$. Booth and Lueker showed that given an interval graph $G$, we can compute the corresponding interval representation in linear time \cite{Booth}. For an interval $I(u)=[l(u),r(u)]$, we say $l(u)$ is the left endpoint and $r(u)$ is the right endpoint of $I(u)$. Without loss of generality, we assume that in the corresponding interval representation of an interval graph, all the endpoints are distinct. An interval graph $G$ is called a \emph{proper interval graph} if, in the corresponding interval representation of $G$, there is no interval that contains some other interval. 

An undirected path graph is the intersection graph of a family of undirected paths of a tree. In \cite{GAVRIL1975237}, Gavril proved that a graph $G$ is an undirected path graph if and only if there exists a tree $\mathcal{T}$ whose vertices are the maximal cliques of $G$ and every maximal clique containing a specific vertex of $G$ forms a path in $\mathcal{T}$. The tree $\mathcal{T}$ is called the \emph{clique tree} of $G$. Next, we propose two linear time algorithms for block graphs and proper interval graphs respectively.

\section{Algorithms}
\subsection{Algorithm in block graphs}\label{sec:Block graph}

In this section, we propose a linear time algorithm to find a minimum liar's ve-dominating set of block graphs. Our algorithm is a labeling-based greedy algorithm. Let $G=(V,E)$ be a block graph with the following labeling on vertices and edges of $G$: $(i)$ $M_V(v)=(t(v),s(v))$, for every vertex $v\in V(G)$, where $t(v)\in \{B, R\}$ and $s(v)$ is a non-negative integer, and $(ii)$ $k(e)\in \{0,1,2\}$, for every edge $e\in E(G)$. A graph $G$ with such a labeling is called a labeled block graph. For any $S\subseteq V$, let $r(S)$ be the set of vertices in $S$ with $t$-label $R$. Next, we define a special dominating set of a labeled block graph $G$ which is essential to design our algorithm. 
\begin{definition}
	Given a labeled block graph $G=(V,E)$, a $M_{LVE}$-dominating set of $G$ is a subset $L\subseteq V(G)$ such that
	\begin{itemize}
		\item[(i)] If $t(v)=R$, then $v\in L$,
		\item[(ii)] for every vertex $v\in V(G)$, $|N_G[v]\cap L|\geq s(v)$,
		\item[(iii)] for every edge $e\in E(G)$, $|N_{G}[e]\cap L|\geq k(e)$, 
		\item[(iv)] for every pair of distinct edges $e\in E(G)$ and $f\in E(G)$, $|(N_G[e]\cup N_G[f])\cap L|\geq \max\{k(e)+k(f)-1,0\}$. 
	\end{itemize}
\end{definition}
The cardinality of the minimum $M_{LVE}$-dominating set of a labeled block graph $G$ is called $M_{LVE}$-domination number of $G$. Note that, if $M_V(v)=(B,0)$ for every vertex $v\in V(G)$ and $k(e)=2$ for every edge $e\in E(G)$, then a minimum $M_{LVE}$-dominating set is a minimum lair's ve-dominating set in $G$. 

\subsubsection{Description of the algorithm}
Now, we give a brief description of our algorithm. Let $G=(V,E)$ be a labeled block graph with the labeling as follows: for every vertex $v\in V(G)$, $M_V(v)=(B,0)$ and for every edge $e\in E(G)$, $k(e)=2$. Also, assume that $\{B_1,B_2,\ldots, B_x\}$ and $\{c_1,c_2,\ldots,c_y\}$ are the blocks and cut vertices of $G$, respectively. We start with an empty $M_{LVE}$-dominating set $L$ of $G$. Next, we construct a rooted cut-tree $T_G$ of $G$ where the root of the cut-tree is a cut vertex of $G$. The set of vertices and edges of $T_G$ are $\{B_1,B_2,\ldots, B_x,c_1,c_2,\ldots,c_y\}$ and $\{B_ic_j:c_j ~is~ a ~cut~vertex ~of ~B_i~ \&~ 1\leq i\leq x ~ \& ~ 1\leq j\leq y \}$, respectively. Observe that, the leaves of the tree $T_G$ are the end blocks of $G$.
Since every block of $G$ is a vertex of $T_G$, we find an ordering $\sigma_{\mathcal{B}}=(B_1,B_2,\ldots,B_x)$ of the blocks in $G$ such that if $i<j$, then $B_i$ appears after $B_j$ in a BFS ordering of the vertices of $T_G$.
If $B_i$ and $B_j$ are two adjacent blocks of $G$ and $B_j$ is the parent of the parent of $B_i$ in $T_G$, then we say $B_i$ is the descendant block of $B_j$. A block of $G$ is called a support block if its' descendant blocks are end blocks of $G$.
Let $B'\in \sigma_{\mathcal{B}}$ be a support block of $G$ and $C=\{c_i:c_i\in V(B')\}$ be the set of cut vertices of $B'$. Suppose that $c_h$ is the cut vertex of $B'$ having a minimum distance from the root of the cut-tree $T_G$. Also, assume that $\{B_{i1},B_{i2},\ldots, B_{ii_q}\}$ is the set of end blocks incident to a cut vertex $c_i$ of $B'$. 

The algorithm is an iterative algorithm where we iterate over the support blocks of $G$. Every cut vertex $c_i\in C\setminus\{c_h\}$ of each support block $B'\in \sigma_{\mathcal{B}}$ is processed in two rounds. In the first round(lines $6$-$18$ of Algorithm \ref{Algo:Blockmain}), we consider the end blocks containing the cut vertex $c_i$ of $B'$. Let $B_{ij}$ be an end block containing the cut vertex $c_{i}$. Let $l_2(B_{ij})$ and $l_1(B_{ij})$ be the number of edges that are only incident to the non-cut vertices of an end block $B_{ij}$ and have $k$-label $2$ and $1$, respectively. We compute $l_2(B_{ij})$, $l_1(B_{ij})$ and $|r(V(B_{ij}))|$. 
To satisfy the conditions $(iii)$ and $(iv)$ of $M_{LVE}$-dominating set for every edge of $B_{ij}$ which is not incident to the cut vertex $c_i$, based on $l_2(B_{ij})$, $l_1(B_{ij})$ and $|r(V(B_{ij}))|$ we relabel some vertices of $V(B_{ij})$ as $R$. Next, we find $s_{max}$ which is the maximum $s$-label of the vertices in $V(B_{ij})\setminus \{c_i\}$ and relabel $\max\{s_{max}-|r(V(B_{ij}))|,0\}$ vertices of $V(B_{ij})$ as $R$.

After the first round, every edge of every end block $B_{ij}$ of $B'$ which is incident to a pair of non-cut vertices of $B_{ij}$ satisfies the conditions $(iii)$ and $(iv)$ of $M_{LVE}$-dominating set. Also, every non-cut vertex of the end blocks of $B'$ satisfies condition $(ii)$ of $M_{LVE}$-dominating set. Thus, only the edges incident to the cut vertex of the end blocks remain to satisfy conditions $(iii)$ and $(iv)$ of $M_{LVE}$-dominating set. In the following round, we select some vertices to satisfy the conditions $(iii)$ and $(iv)$ of $M_{LVE}$-dominating set for the edges of end blocks of $B'$ incident to the cut vertex $c_i$.

\begin{algorithm}
	\scriptsize
	\textbf{Input:} A labelled block graph $G=(V,E)$ with $M_V(v)=(B,0)$ for every $v\in V$ and $k(e)=2$ for every $e\in E$.\\
	\textbf{Output:} A minimum liar's ve-dominating set $L$ of $G$.
	\caption{LVED-Block$(G)$}
	\label{Algo:Blockmain}
	\begin{algorithmic}[1]
		\State $L=\phi$;
		\State Construct a rooted cut tree $T_G$ with a cut vertex as root;
		\State $\sigma_{\mathcal{B}}=(B_1,B_2,\ldots, B_s)$ such that $B_s$ is a block at level $1$ of $T_G$;
		\For{(every support block $B'$ in $\sigma_{\mathcal{B}}$)}
		\For{(every cut vertex $c_{i}\in C\setminus\{c_h\}$ of $B'$)}
		\For{(every end block $B_{ij}$ incident to $c_i$)}
		\State Compute $l_2(B_{ij}),l_1(B_{ij}),r(V(B_{ij}))$;
		\If{($l_2(B_{ij})\geq 2$)}
		\State Relabel $\max\{3-|r(V(B_{ij}))|,0\}$ vertices of $V(B_{ij})$ as $R$ with priority $c_i$;
		\State $|r(V(B_{ij}))|=|r(V(B_{ij}))|+\max\{3-|r(V(B_{ij}))|,0\}$;
		\ElsIf{($l_2(B_{ij})=1$)}
		\State Relabel $\max\{2-|r(V(B_{ij}))|,0\}$ vertices of $V(B_{ij})$ as $R$ with priority $c_i$;
		\State $|r(V(B_{ij}))|=|r(V(B_{ij}))|+\max\{2-|r(V(B_{ij}))|,0\}$;
		\ElsIf{($l_2(B_{ij})=0 ~ \& ~l_1(B_{ij})\geq 1$)}
		\State Relabel $\max\{1-|r(V(B_{ij}))|,0\}$ vertices of $V(B_{ij})$ as $R$ with priority $c_i$;
		\State $|r(V(B_{ij}))|=|r(V(B_{ij}))|+\max\{1-|r(V(B_{ij}))|,0\}$;
		\EndIf 
		\State $s_{max}=\max\{s(v):v\in V(B_{ij})\setminus \{c_i\}\}$;	
		\State Relabel $\max\{s_{max}-|r(V(B_{ij}))|,0\}$ vertices of $V(B_{ij})$ as $R$ with priority $c_i$;
		\EndFor
		\EndFor
		
		\For{(every cut vertex $c_i\in C\setminus\{c_h\}$ of $B'$)}
		\State Compute $l_2(c_i), r(N_G[c_i]), max_1(c_i),max_2(c_i), w(c_i)$;
		\If{($w(c_i)>0$)}
		\State Relabel $w(c_i)$ vertices of $N_G[c_{i}]$ as $R$ with priority the vertices of $V(B')$;
		\EndIf  
		\EndFor
		\If{(there exists $c_a,c_b\in C\setminus\{c_h\}$ with $l_2(c_{a})=l_2(c_{b})=1$ \& $|r(N_G(c_a)\setminus V(B'))|=|r(N_G(c_b)\setminus V(B'))|=0$)}
		\State Relabel $\max\{(3-|r(V(B'))|),0\}$ vertices of $V(B')$ as $R$ with priority $c_h$;
		\EndIf
		
		\State $q=$ index of a cut vertex in $C\setminus\{c_h\}$ such that $|r(N_G[c_q])|=\min\{|r(N_G[c_i])| : c_i\in C\setminus\{c_h\}~and~l_2(c_i)=1\}$;
		\If{($V(B')\setminus C\neq \emptyset$)}
		\State Relabel $\max\{3-|r(N_G[c_q])|,0\}$ vertices of $V(B')$ as $R$ with priority $c_h$;
		\ElsIf{(there exists a vertex $c_i\in C\setminus\{c_h,c_q\}$ such that $r(N_G(c_i)\setminus V(B'))=\emptyset$)}
		\State $s(c_i)=\max\{s(c_i), 3-|r(N_G[c_q]|\}$;
		\EndIf
		\State $s(c_h)=\max\{s(c_h),3-|r(N_G[c_q]|\}$;	
		\State $L=L\bigcup\limits_{c_i\in B'} r(N_G(c_i)\setminus V(B')) $;
		\State Delete all the end blocks adjacent to $B'$;
		\For{(every edge $e$ incident to cut vertices $c_i$ of $B'$)} 			
		\State $k(e)=\max\{k(e)-|r(N_G(c_i)\setminus V(B'))|,0\}$;
		\EndFor
		\For{(every cut vertex $c_i$ of $B'$)} 			
		\State $s(c_i)=\max\{s(c_i)-|r(N_G(c_i)\setminus V(B'))|,0\}$;
		\EndFor
		\EndFor
		\State $D=$StarBlock$(G)$;
		\State $L=L\cup D$;
		\State \Return $L$;
	\end{algorithmic}
\end{algorithm}

In the second round(lines $19$-$22$ of Algorithm \ref{Algo:Blockmain}), we again process the cut vertices $c_i\in C\setminus\{c_h\}$ of $B'$. Let $max_1(c_i)$ and $max_2(c_i)$ be the maximum and second maximum $k$-label of the edges of the end blocks of $B'$ incident to $c_i$. In this round, we first find $max_1(c_i)$ and $max_2(c_i)$ and we relabel $w(c_i)$ vertices of $B'$ where $w(c_i)=\max\{max_1(c_i)-|r(N_G[c_i])|, \max\{max_1(c_i)+max_2(c_i)-1,0\}-|r(N_G[c_i])|\}$. In other words, $w(c_i)$ is the number of vertices that need to be selected to satisfy condition $(iii)$ for every edge of the end blocks incident to a cut vertex $c_i$ and condition $(iv)$ for every pair of edges of the end blocks incident to the cut vertex $c_i$. After the second round, every edge of the end blocks incident to the cut vertices $C\setminus\{c_h\}$ of $B'$ satisfies the condition $(iii)$ and every pair of edges of the same end block satisfies the condition $(iv)$ of $M_{LVE}$-dominating set.

\begin{algorithm}
	\scriptsize
	\textbf{Input:} A labelled star block graph $G=(V,E)$ with cut vertex $c'$.\\
	\textbf{Output:} A minimum $M_{LVE}$-dominating set $D$ of $G$.
	\caption{StarBlock$(G)$}
	\label{Algo:StarBlock}
	\begin{algorithmic}[1]
		\For{(every end block $B_{ij}$ of $G$)}
		\If{($l_2(B_{ij})\geq 2$)}
		\State Relabel $3-|r(V(B_{ij}))|$ vertices of $V(B_{ij})$ as $R$ with priority $c'$;
		\ElsIf{($l_2(B_{ij})=1$)}
		\State Relabel $2-|r(V(B_{ij}))|$ vertices of $V(B_{ij})$ as $R$ with priority $c'$;
		\ElsIf{($l_2(B_{ij})=0 ~ \& ~ l_1(B_{ij})=1$)}
		\State Relabel $1-|r(V(B_{ij}))|$ vertices of $V(B_{ij})$ as $R$ with priority $c'$;
		\EndIf
		\State Relabel $\max\{s_{max}-|r(V(B_{ij}))|,0\}$ vertices of $V(B_{ij})$ as $R$ with priority $c'$; 	
		\EndFor
		\State Compute $l_2(c'), r(N_G[c']), max_1(c'),max_2(c'), w(c')$;
		\If{($w(c')>0$)}
		\State Relabel $w(c')$ vertices of $N_G[c']$ as $R$ with priority $c'$;
		\EndIf 
		\State Relabel $\max\{s(c')-|r(N_G[c'])|,0\}$ vertices of $N_G[c']$ as $R$;
		\State $D=D\cup r(c')$;
		\State \Return $D$ ;
	\end{algorithmic}
\end{algorithm}

Let $l_2(c_i)$ be the number of edges of the end blocks of $B'$ incident to $c_i$ with $k$-label $2$. Observe that if $|r(N_G[c_i]\setminus V(B'))|>0$ and $l_2(c_i)=1$, the condition $(iv)$ holds for every edge of the end blocks of $B'$ incident to $c_i$ with $k$-label $2$. Thus, we check if there is any pair of cut vertices with $l_2(c_i)=1$ and $|r(N_G[c_i]\setminus V(B'))|=0$. In that case, we relabel some vertices of $B'$ as $R$, if needed. After that, every pair of edges of end blocks of $B'$ incident to the cut vertices of $B'$ satisfies condition $(iv)$. This implies that the only pair of edges for which the condition $(iv)$ does not hold is an edge of the end block incident to the cut vertices of $B'$ with $k$-label $2$ and an edge of $B'$. Due to our choice of vertices in the second round of processing the cut vertices of $B'$, the edges incident to the cut vertices $c_i\in C\setminus\{c_h\}$ of $B'$ with $l_2(c_i)=2$ satisfy condition $(iv)$ of $M_{LVE}$-dominating set. Thus, we only need to check the cut vertices of $B'$ with $l_2(c_i)=1$.

Next, we find the cut vertex $c_q$ of $B'$ with $l_2(c_q)=1$ which has the minimum number of vertices with $t$-label $R$ in its neighbourhood. Now, if the support block $B'$ contains a non-cut vertex, then to satisfy the condition $(iv)$ we update the $t$-label of $\max\{3-|r(N_G[c_q])|,0\}$ many vertices of $B'$ as $R$. Also, if there is a cut vertex of $B'$ with no vertex of $t$-label $R$ in the end blocks containing it, then to satisfy condition $(iv)$ for the pair of edges(one incident to that cut vertex and $c_q$ and another is the edge with $k$-label $2$ incident to the vertex $c_q$) we update the $s$-label of that cut vertex by $s(c_i)=\max\{s(c_i), 3-|r(N_G[c_q]|\}$. Now, we update the $s$-label of the cut vertex $c_h$ to ensure that the pair of edges $c_hc_q$ and the edge with $k$-label $2$ incident to $c_q$ satisfy the condition $(iv)$ of $M_{LVE}$-dominating set. After that, we include all the vertices(except the cut vertices) of the end blocks of $B'$ with $t$-label $R$ in the set $L$.

Next, we delete all the end blocks incident to the cut vertices of $C\setminus\{c_h\}$ and update the $k$-label and $s$-label of the edges incident to the cut vertices and the cut vertices of $B'$, respectively. After that, we proceed to the next iteration. Repeating these steps a finite number of times yields a labeled star block graph $G'$. Now, we call the algorithm StarBlock($G'$)(Algorithm \ref{Algo:StarBlock}). StarBlock() algorithm similarly deals with the end blocks like Algorithm \ref{Algo:Blockmain}. The only difference is that after processing all the end blocks, in Algorithm \ref{Algo:StarBlock}, we check if the cut vertex $c'$ satisfies the condition $(ii)$ of $M_{LVE}$-domination or not and based on that we relabel some vertices of $N_G[c']$ as $R$(line $12$ of Algorithm \ref{Algo:StarBlock}). the algorithm StarBlock($G'$) returns a minimum $M_{LVE}$-dominating set of $G'$. We include those vertices in the set $L$ which is the minimum liar's ve-dominating set of $G$. The details of the algorithm are in Algorithm \ref{Algo:Blockmain}.

\subsubsection{Correctness}
Now, we prove the correctness of Algorithm \ref{Algo:Blockmain} and Algorithm \ref{Algo:StarBlock}. Let $\gamma(G)$ be the $M_{LVE}$-domination number of a block graph $G$. First, we prove the correctness of Algorithm \ref{Algo:Blockmain}. 

\begin{lem}\label{lem:1if}
	Let $B'$ be a support block in $\sigma_{\mathcal{B}}$, $c_i$ be a cut vertex of $B'$ processed in the first round, and $B_{ij}$ be an end block containing the vertex $c_i$. 
	
	\noindent(a) If $l_2(B_{ij})\geq 2$, then $\gamma(G)=\gamma(G')$, where $G'$ is obtained by relabelling $\max\{(3-|r(V(B_{ij}))|),0\}$ vertices of $B_{ij}$ with priority $c_i$ as $R$ and every other label remains the same.
	
	\noindent(b) If $l_2(B_{ij})= 1$, then $\gamma(G)=\gamma(G')$, where $G'$ is obtained by relabelling $\max\{(2-|r(V(B_{ij}))|),0\}$ vertices of $B_{ij}$ with priority $c_i$ as $R$ and every other label remains the same. 
	
	\noindent(c) If $l_2(B_{ij})=0$ and $l_1(B_{ij})\geq 1$, then $\gamma(G)=\gamma(G')$, where $G'$ is obtained by relabelling $\max\{(1-|r(V(B_{ij})|)),0\}$ vertices of $B_{ij}$ with priority $c_i$ as $R$ and every other label remains the same.
	
	\noindent(d) $\gamma(G)=\gamma(G')$, where $G'$ is obtained by relabelling $s_{max}-|r(V(B_{ij}))|$ vertices of $B_{ij}$ as $R$ with priority $c_i$ and every other label remains the same.
\end{lem}

\begin{proof}
	\noindent(a) Let $L$ be a minimum $M_{LVE}$-dominating set of $G$. Since $l_2(B_{ij})\geq 2$, there exists at least two edges of $B_{ij}$ which are not incident to $c_i$, say $e_p$ and $e_q$, with $k(e_p)=k(e_q)=2$. Now, by the definition of $L$, we have $|(N_G[e_p]\cup N_G[e_q])\cap L|\geq \max\{k(e_p)+k(e_q)-1,0\}=3$. Observe that $N_G[e_p]$ and $N_G[e_q]$ are the subset of vertices of $B_{ij}$. Since every vertex $v\in B_{ij}$ with $t(v)=R$ is in $L$, we have $r(V(B_{ij}))\subset L$. Therefore, $L$ must contain at least $\max\{(3-|r(V(B_{ij}))|),0\}$ vertices of $B_{ij}$ having $t$-label $B$ in $G$. Let $L'$ be the set obtained from $L$ by replacing $\max\{(3-|r(V(B_{ij}))|),0\}$ vertices of $B_{ij}$ having $t$-label $B$ in $G$ with the vertices of $B_{ij}$ having $t'$-label $R$ in $G'$. Since every other label remains the same, it is easy to observe that $L'$ is a $M_{LVE}$-dominating set of $G'$. Hence, we have $\gamma(G')\leq \gamma(G)$. 
	
	Conversely, let $L'$ be a minimum $M_{LVE}$-dominating set of $G'$. Since every edge of $G'$ has same $k$-label as $G$ and every vertex of $G'$ has the same $M_V$-label as $G$ except for $\max\{(3-|r(V(B_{ij}))|),0\}$ many vertices of $B_{ij}$. Therefore, $L'$ is a $M_{LVE}$-dominating set of $G$. Hence, we have 
	$\gamma(G)\leq \gamma(G')$. Thus, $\gamma(G)=\gamma(G')$. \qed
	
	\noindent(b) Let $L$ be a minimum $M_{LVE}$-dominating set of $G$. Since $l_2(B_{ij})=1$, there exists exactly one edge of $B_{ij}$ which is not incident to $c_i$, say $e_p$, with $k(e_p)=2$. Now, by the definition of $L$, we have $|N_G[e_p]\cap L|\geq k(e_p)=2$. Observe that $N_G[e_p]$ is the set of vertices of $B_{ij}$. Since every vertex $v\in B_{ij}$ with $t(v)=R$ is in $L$, we have $r(V(B_{ij}))\subset L$. Therefore, $L$ must contain at least $\max\{(2-|r(V(B_{ij}))|),0\}$ vertices of $B_{ij}$ having $t$-label $B$ in $G$. Let $L'$ be the set obtained from $L$ by replacing $\max\{(2-|r(V(B_{ij}))|),0\}$ vertices of $B_{ij}$ having $t$-label $B$ in $G$ with the vertices of $B_{ij}$ having $t'$-label $R$ in $G'$. Since every other label remains the same, it is easy to observe that $L'$ is a $M_{LVE}$-dominating set of $G'$. Hence, we have $\gamma(G')\leq \gamma(G)$. 
	
	Conversely, let $L'$ be a minimum $M_{LVE}$-dominating set of $G'$. Since every edge of $G'$ has same $k$-label as $G$ and every vertex of $G'$ has the same $M_V$-label as $G$ except for $\max\{(2-|r(V(B_{ij}))|),0\}$ many vertices of $B_{ij}$. Therefore, $L'$ is a $M_{LVE}$-dominating set of $G$. Hence, we have $\gamma(G)\leq \gamma(G')$. Thus, we have $\gamma(G)=\gamma(G')$.\qed
	
	\noindent(c) Let $L$ be a minimum $M_{LVE}$-dominating set of $G$. Since $l_2(B_{ij})=0$ and $l_1(B_{ij})\geq 1$, the maximum $k$-label of the edges of $B_{ij}$ which is not incident to $c_i$ is $1$. Let $e_p$ be such an edge with $k(e_p)=1$. Now, by the definition of $L$, we have $|N_G[e_p]\cap L|\geq k(e_p)=1$. Observe that $N_G[e_p]$ is the set of vertices of $B_{ij}$. Since every vertex $v\in B_{ij}$ with $t(v)=R$ is in $L$, we have $r(V(B_{ij}))\subset L$. Therefore, $L$ must contain at least $\max\{(1-|r(V(B_{ij}))|),0\}$ vertices of $B_{ij}$ having $t$-label $B$ in $G$. Let $L'$ be the set obtained from $L$ by replacing $\max\{(1-|r(V(B_{ij}))|),0\}$ vertices of $B_{ij}$ having $t$-label $B$ in $G$ with the vertices of $B_{ij}$ having $t'$-label $R$ in $G'$. Since every other label remains the same, it is easy to observe that $L'$ is a $M_{LVE}$-dominating set of $G'$. Hence, we have $\gamma(G')\leq \gamma(G)$. 
	
	Conversely, let $L'$ be a minimum $M_{LVE}$-dominating set of $G'$. Since every edge of $G'$ has same $k$-label as $G$ and every vertex of $G'$ has the same $M_V$-label as $G$ except for $\max\{(1-|r(V(B_{ij}))|),0\}$ many vertices of $B_{ij}$. Therefore, $L'$ is a $M_{LVE}$-dominating set of $G$. Hence, we have $\gamma(G)\leq \gamma(G')$. Thus, we have $\gamma(G)=\gamma(G')$.\qed
	
	\noindent(d) Let $L$ be a $M_{LVE}$-dominating set of $G$. Thus, for every vertex $v\in V$, we have $|N_G[v]\cap L|\geq s(v)$. Let $s_{max}=\max\{s(v ):v\in V(B_{ij})\setminus\{c_i\}\}$. Hence, $|N_G[v]\cap L|\geq s_{max}$ for any $v\in V(B_{ij})$. Now, for any vertex $v\in V(B_{ij})\setminus\{c_i\}$, $N_G[v]$ is a subset of $V(B_{ij})$. Since every vertex with $t$-label $R$ is in $L$, we have $r(V(B_{ij}))\subseteq L$. This implies that $L$ must contain at least $\max\{(s_{max}-|r(V(B_{ij}))|),0\}$ vertices of $V(B_{ij})$ having $t$-label $B$. Let $L'$ be the set obtained from $L$ by replacing the $\max\{(s_{max}-|r(V(B_{ij}))|),0\}$ vertices of $V(B_{ij})$ having $t$-label $B$ in $G$ with the $\max\{(s_{max}-|r(V(B_{ij}))|),0\}$ vertices of $V(B_{ij})$ having $t'$-label $R$ in $G'$. Since every other label remains the same, it is easy to observe that $L'$ is a $M_{LVE}$-dominating set of $G'$. Hence, we have $\gamma(G')\leq \gamma(G)$. 
	
	Conversely, let $L'$ be a minimum $M_{LVE}$-dominating set of $G'$. Since every edge of $G'$ has same $k$-label as $G$ and every vertex of $G'$ has the same $M_V$-label as $G$ except for $\max\{(s_{max}-|r(V(B_{ij}))|),0\}$ many vertices of $B_{ij}$. Therefore, $L'$ is a $M_{LVE}$-dominating set of $G$. Hence, we have $\gamma(G)\leq \gamma(G')$. Thus, we have $\gamma(G)=\gamma(G')$.
	%
\end{proof}

From Lemma \ref{lem:1if}, it follows that after the first round the labeling of $G$ is such that for any edge $e$ of any descendant block $B_{ij}$ of $B'$ which is not incident to the cut vertex $c_i$ of $B_{ij}$, $k(e)\leq |r(N_G[c_i])|$, for every pair of edges $(e,e')$ of every end block $B_{ij}$ of $B'$ which is not incident to the cut vertex $c_i$ of $B_{ij}$, $\max\{k(e)+k(e')-1,0\}\leq |r(N_G[c_i])|$ and for every non-cut vertex $v$ of every end block $B_{ij}$, $s(v)<|r(N_G[c_i])|$ where $c_i$ is the cut vertex of $B_{ij}$. Considering these assumptions, we now prove the correctness of the algorithm for the second round of processing of the cut vertices of the support block $B'$. 


\begin{lem}\label{lem:2if}
	Let $B'$ be a support block in $\sigma_{\mathcal{B}}$, $c_i$ be a cut vertex of $B'$ processed at the second round and $B_{ij}$ be any end block containing the vertex $c_i$. Also, assume that labels of $G$ are such that for every edge $e$ of $B_{ij}$ which is not incident to $c_i$, $k(e)\leq |r(N_G[c_i])|$, for every such pair of edges $(e,e')$ of $B_{ij}$, $\max\{k(e)+k(e')-1,0\}\leq |r(N_G[c_i])|$ and for every vertex $v\in V(B_{ij})\setminus\{c_i\}$, $s(v)<|r(N_G[c_i])|$. Let $max_1(c_i)$ and $max_2(c_2)$ be the maximum and second maximum $k$-label of the edges between $c_i$ and the end blocks incident to $c_i$. If $w(c_i)>0$, where $w(c_i)=\max\{max_1(c_i)-|r(N_G[c_i])|, \max\{max_1(c_i)+max_2(c_i)-1,0\}-|r(N_G[c_i])|\}$, then $\gamma(G)=\gamma(G')$, where $G'$ is obtained by relabelling $w(c_i)$ vertices of $N_G[c_i]$ with priority to the vertices of $V(B')$ as $R$ and every other label remains the same.
\end{lem}

\begin{proof}
	Let $L$ be a minimum $M_{LVE}$-dominating set of $G$. Observe that, the neighbourhood of any edge $e$ of the descendant blocks of $B'$ incident to $c_i$, $N_G[e]$ is a subset of $N_G[c_i]$. Hence, by the definition of $M_{LVE}$-dominating set, $L$ contains at least $\max\{max_1(c_i), max_1(c_i)+max_2(c_i)-1\}$ many vertices from $N_G[c_i]$. Since $L$ must contain every vertex with $t$-label $R$,  $r(N_G[c_i])$ is a subset of $L$. Now, $w(c_i)=\max\{
	max_1(c_i)-|r(N_G[c_i])|, \max\{max_1(c_i)+max_2(c_i)-1,0\}-|r(N_G[c_i])|\}$ and $w(c_i)>0$. Therefore, $L$ contains $w(c_i)$ vertices of $N_G[c_i]$ with $t$-label $B$. Let $L'$ be the set obtained from $L$ by replacing the $w(c_i)$ vertices of $N_G[c_i]$ with $t$-label $B$ in $G$ with the vertices of $N_G[c_i]$ having $t'$-label as $R$ in $G'$. Since every other label remains the same, it is easy to verify that $L'$ is a $M_{LVE}$-dominating set of $G'$. Hence, we have $\gamma(G')\leq \gamma(G)$. 
	
	Conversely, let $L'$ be a minimum $M_{LVE}$-dominating set of $G'$. Since every edge of $G'$ has same $k$-label as $G$ and every vertex of $G'$ has the same $M_V$-label as $G$ except for $w(c_i)$ vertices of $N_G[c_i]$. Therefore, $L'$ is a $M_{LVE}$-dominating set of $G$. Hence, we have $\gamma(G)\leq \gamma(G')$. Thus, $\gamma(G)=\gamma(G')$.
\end{proof}

After the second round of processing of the cut vertices of $B'$, the labeling of $G$ is as follows: 
for every end block $B_{ij}$ containing any cut vertex $c_i\in C\setminus\{c_h\}$ of $B'$, $k(e)\leq |r(N_G[c_i])|$ for every edge $e$ of $B_{ij}$, $\max\{k(e)+k(e')-1,0\}\leq |r(N_G[c_i])|$ for every pair of edges $(e,e')$ of $B_{ij}$, and $s(v)\leq |r(N_G[c_i])|$ for every non-cut vertex $v$ of $B_{ij}$. Thus, the condition $(iii)$ of $M_{LVE}$-dominating set holds for every edge of the end blocks of $B'$.
Now, we show that condition $(iv)$ of $M_{LVE}$-domination holds for every pair of edges of descendant blocks of $B'$. Note that, from Lemma \ref{lem:1if} and Lemma \ref{lem:2if}, the condition $(iv)$ holds for the pair of edges of the end blocks containing a cut vertex $c_i$.  Therefore, it is sufficient to show that condition $(iv)$ holds for the pair of edges of the end blocks incident to different cut vertices of the support block $B'$. We prove this in the next lemma.  

\begin{lem}\label{lem:3if}
	Let $B'$ be a support block in $\sigma_{\mathcal{B}}$, $c_i\in C\setminus\{c_h\}$ is a cut vertex of $B'$ and $B_{ij}$ is an end block containing the vertex $c_i$. Also, assume that the labels of $G$ are as follows:  $k(e)\leq |r(N_G[c_i])|$ for every edge $e$ of $B_{ij}$, $\max\{k(e)+k(e')-1,0\}\leq |r(N_G[c_i])|$ for every pair of edges $(e,e')$ of $B_{ij}$, and $s(v)\leq |r(N_G[c_i])|$ for every non-cut vertex $v$ of $B_{ij}$. If there is a pair of cut vertices $c_a,c_b\in C\setminus\{c_h\}$ of $B'$ with $l_2(c_a)=l_2(c_b)=1$ and $|r(N_G(c_a)\setminus V(B'))|=|r(N_G(c_b)\setminus V(B'))|=0$, then $\gamma(G)=\gamma(G')$, where $G'$ is obtained by relabelling $\max\{(3-|r(V(B'))|),0\}$ vertices of $B'$ with priority $c_h$ as $R$ and every other label remains the same. 
\end{lem}

\begin{proof}
	%
	Let $L$ be a minimum $M_{LVE}$-dominating set of $G$. Since $l_2(c_a)=l_2(c_b)=2$ for the cut vertices $c_a$ and $c_b$, there are two edges, say $e_a$ and $e_b$, of two end blocks with $k$-label $2$ and incident to $c_a$ and $c_b$, respectively. Therefore, $L$ contains at least $\max\{k(e_a)+k(e_b)-1,0\}=3$ vertices of $N_G[c_a]\cup N_G[c_b]$. Now, every vertex with $t$-label $R$ is in $L$. Therefore, $L$ contains all the vertices of $r(V(B'))$. Since $r(N_G[c_a]\setminus V(B'))=r(N_G[c_b]\setminus V(B'))=\emptyset$, the only vertices of $N_G[c_a]\cup N_G[c_b]$ with $t$-label $R$ is in $V(B')$. This implies that $L$ contains at least $\max\{(3-|r(V(B'))|),0\}$ vertices of $N_G[c_a]\cup N_G[c_b]$ with $t$-label $B$. Let $L'$ be the set obtained from $L$ by replacing the $\max\{(3-|r(V(B'))|),0\}$ vertices of $N_G[c_a]\cup N_G[c_b]$ having $t$-label $B$ in $G$ with the $\max\{(3-|r(V(B'))|),0\}$ vertices of $V(B')$ having $t'$-label $R$ in $G'$. Since every other label remains the same, it is easy to observe that $L'$ is a $M_{LVE}$-dominating set of $G'$. Hence, we have $\gamma(G')\leq \gamma(G)$. 
	
	Conversely, let $L'$ be a minimum $M_{LVE}$-dominating set of $G'$. Since every edge of $G'$ has same $k$-label as $G$ and every vertex of $G'$ has the same $M_V$-label as $G$ except for $\max\{(3-|r(V(B'))|),0\}$ many vertices of $V(B')$. Therefore, $L'$ is a $M_{LVE}$-dominating set of $G$. Hence, we have $\gamma(G)\leq \gamma(G')$. Thus, we have $\gamma(G)=\gamma(G')$.	
\end{proof}

From Lemma \ref{lem:3if}, it follows that the condition $(iii)$ of $M_{LVE}$-domination holds for every edge of the end blocks of $B'$ and condition $(iv)$ of $M_{LVE}$-domination holds for every pair of edges of the end blocks of $B'$. Let $B_{ij}$ be any end block of $B'$ containing the cut vertex $c_i$ of $B'$. The labeling of $G$ is as follows: $k(e)\leq |r(N_G[c_i])|$ for every edge $e$ of $B_{ij}$, $\max\{k(e)+k(e')-1,0\}\leq |r(N_G[c_i])|$ for every pair of edges $e,e'$ of $B_{ij}$, $s(v)\leq |r(N_G[c_i])|$ for every non-cut vertex $v$ of $B_{ij}$, and $\max\{k(e)+k(e')-1,0\}\leq |r(N_G[c_a]\cup N_G[c_b])|$ for every pair of edges $(e,e')$ such that $e$ and $e'$ are the edges of the end blocks $B_{aj}$ and $B_{bj'}$, where $B_{aj}$ and $B_{bj'}$ are end blocks containing the cut vertices $c_a$ and $c_b$, respectively. Now, we show that the condition $(iv)$ holds for every pair of edges where one edge is from the end block of $B'$ and another is incident to a non-cut vertex of $B'$.

\begin{lem}\label{lem:4}
	Let $B'$ be a support block in $\sigma_{\mathcal{B}}$, $c_i\in C\setminus\{c_h\}$ is a cut vertex of $B'$ and $B_{ij}$ is an end block containing the vertex $c_i$. Also, assume that the labels of $G$ are as follows: $k(e)\leq |r(N_G[c_i])|$ for every edge $e$ of $B_{ij}$, $\max\{k(e)+k(e')-1,0\}\leq |r(N_G[c_i])|$ for every pair of edges $e,e'$ of $B_{ij}$, $s(v)\leq |r(N_G[c_i])|$ for every non-cut vertex $v$ of $B_{ij}$, and $\max\{k(e)+k(e')-1,0\}\leq |r(N_G[c_a]\cup N_G[c_b])|$ for every pair of edges $(e,e')$ such that $e$ and $e'$ are the edges of the end blocks $B_{aj}$ and $B_{bj'}$, where $B_{aj}$ and $B_{bj'}$ are end blocks containing the cut vertices $c_a$ and $c_b$, respectively. Let $q$ be the index of a cut vertex in $C\setminus\{c_h\}$ such that $|r(N_G[c_q])|$ is minimum of all $|r(N_G[c_i])|$ corresponding to every cut vertices $c_i\in C\setminus\{c_h\}$ having $l_2(c_i)=1$. If $B'$ has a non-cut vertex, then, $\gamma(G)=\gamma(G')$, where $G'$ is obtained by relabelling $\max\{3-|r(N_G[c_q])|,0\}$ vertices of $B'$ as $R$ with priority $c_h$ and every other label remains the same.
\end{lem}

\begin{proof}
	Let $L$ be a minimum $M_{LVE}$-dominating set of $G$. Since $l_2(c_q)=1$, there is only one edge incident to $c_q$ with $k$-label $2$, say $e_q=c_qv_x$, where $v_x$ is a vertex of an end block of $B'$. Now, for every pair of edges having $k$-label $2$ and incident to $c_q$, $L$ contains at least $3$ vertices to satisfy the condition $(iv)$ of $M_{LVE}$-dominating set. Suppose that, $v'$ is a non-cut vertex of $B'$. Now, consider the edge $e'=c_qv'$. Clearly, $k$-label of $e'$ is $2$, and hence, $L$ contains at least $3$ vertices from $N_G[c_q]\cup N_G[v']$. Since $v'$ is a non-cut vertex of $B'$ and $c_q$ is a cut vertex of $B'$, we have $N_G[v']\subset N_G[c_q]$. Thus, $L$ contains at least $3$ vertices from $N_G[c_q]$. Now, every vertex with $t$-label $R$ is in $L$. Therefore, $L$ contains all the vertices of $r(N_G[c_q])$. This implies that $L$ contains at least $\max\{(3-|r(N_G[c_q])|),0\}$ vertices of $N_G[c_q]$ with $t$-label $B$. Let $L'$ be the set obtained from $L$ by replacing the $\max\{(3-|r(N_G[c_q])|),0\}$ vertices of $N_G[c_q]$ having $t$-label $B$ in $G$ with the $\max\{(3-|r(N_G[c_q])|),0\}$ vertices of $V(B')$ having $t'$-label $R$ in $G'$. Since every other label remains the same, it is easy to observe that $L'$ is a $M_{LVE}$-dominating set of $G'$. Hence, we have $\gamma(G')\leq \gamma(G)$.

	
	Conversely, let $L'$ be a minimum $M_{LVE}$-dominating set of $G'$. Since every edge of $G'$ has same $k$-label as $G$ and every vertex of $G'$ has the same $M_V$-label as $G$ except for $\max\{(3-|r(N_G[c_q])|),0\}$ many vertices of $V(B')$. Therefore, $L'$ is a $M_{LVE}$-dominating set of $G$. Hence, we have $\gamma(G)\leq \gamma(G')$. Thus, we have $\gamma(G)=\gamma(G')$.
	%
\end{proof}

In the following Lemma, we show that condition $(iv)$ of $M_{LVE}$-domination holds for every pair of edges where one edge is from an end block of $B'$ and another is incident to a cut vertex of $B'$ with no vertices of $t$-label $R$ in the end blocks containing it.

\begin{lem}\label{lem:5}
	Let $B'$ be a support block in $\sigma_{\mathcal{B}}$, $c_i\in C\setminus\{c_h\}$ is a cut vertex of $B'$ and $B_{ij}$ is an end block containing the vertex $c_i$. Also, assume that the labels of $G$ are as follows: $k(e)\leq |r(N_G[c_i])|$ for every edge $e$ of $B_{ij}$, $\max\{k(e)+k(e')-1,0\}\leq |r(N_G[c_i])|$ for every pair of edges $e,e'$ of $B_{ij}$, $s(v)\leq |r(N_G[c_i])|$ for every non-cut vertex $v$ of $B_{ij}$, and $\max\{k(e)+k(e')-1,0\}\leq |r(N_G[c_a]\cup N_G[c_b])|$ for every pair of edges $(e,e')$ such that $e$ and $e'$ are the edges of the end blocks $B_{aj}$ and $B_{bj'}$, where $B_{aj}$ and $B_{bj'}$ are end blocks containing the cut vertices $c_a$ and $c_b$, respectively. Let $q$ be the index of a cut vertex in $C\setminus\{c_h\}$ such that $|r(N_G[c_q])|$ is minimum of all $|r(N_G[c_i])|$ corresponding to every cut vertices $c_i\in C\setminus\{c_h\}$ having $l_2(c_i)=1$. If $B'$ has a cut vertex $c_i\in C\setminus\{c_h,c_q\}$ with $r(N_G(c_i)\setminus V(B'))=\emptyset$, then $\gamma(G)=\gamma(G')$, where $G'$ is obtained by relabelling $s(c_i)$ as $s'(c_i)=\max\{s(c_i),3-|r(N_G[c_q])|\}$ and every other label remains the same.
\end{lem}

\begin{proof}
	Let $L$ be a minimum $M_{LVE}$-dominating set of $G$. Since $l_2(c_q)=1$, there is only one edge, say $e_q=c_qv_x$, where $v_x$ is a vertex of an end block incident to $c_q$, with $k$-label $2$. Now, for every pair of edges having $k$-label $2$ and incident to $c_q$, $L$ contains at least $3$ vertices to satisfy the condition $(iv)$ of $M_{LVE}$-dominating set. Consider the edge $e'=c_qc_i$. Clearly, $k$-label of $e'$ is $2$ and hence $L$ contains at least $3$ vertices from $N_G[c_q]\cup N_G[c_i]$. Since $L$ contains every vertex with $t$-label $R$, $r(N_G[c_q]\cup N_G[c_i])$ is contained in $L$. Therefore, $L$ contains at least $3-(|r(N_G[c_q]|)|)$ vertices of $N_G[c_q]\cup N_G[c_h]$ having $t$-label $B$ because $r(N_G(c_i)\setminus V(B'))=\emptyset$ and $N_G[c_i]\cap N_G[c_q]=V(B')$. Also, $L$ contains at least $s(v)$ vertices from $N_G[v]$ for every vertex $v\in V(G)$. Hence, there are $s(c_q)$ and $s(c_i)$ vertices in $L$ corresponding to the vertices $c_q$ and $c_i$, respectively. If $s(c_i)\geq 3-|r(N_G[c_i])|$, then $s(c_i)=s'(c_i)$ and thus, $L$ is a $M_{LVE}$-dominating set of $G'$. Otherwise, let $L'$ be the set obtained from $L$ by replacing the $3-(|r(N_G[c_q]|))$ vertices of $N_G[c_q]\cup N_G[c_i]$ having $t$-label $B$ with the $3-(|r(N_G(c_q))|)$ vertices of $N_G[c_i]\cap V(B')$ having $t$-label $B$. Hence, $L'$ contains at least $\max\{s(c_i),3-(|r(N_G[c_q])|)$ many vertices for the vertex $c_i$. Therefore, it is easy to observe that $L'$ is a $M_{LVE}$-dominating set of $G'$. Hence, we have $\gamma(G')\leq \gamma(G)$. 
	
	Conversely, let $L'$ be a minimum $M_{LVE}$-dominating set of $G'$. Now, every edge of $G'$ has same $k$-label as $G$ and every vertex of $G'$ has the same $M_V$-label as $G$ except for the vertex $c_i$. We show that $L'$ is a $M_{LVE}$-dominating set of $G$. Since $L'$ contains at least $s'(c_i)=\max\{s(c_i),3-(|r(N_{G}[c_q])|)\}$ vertices in $N_{G'}[c_i]$, it is straightforward that $L'$ contains at least $s(c_i)$ vertices in $N_G[c_i]$. Since every other label is the same, it is easy to see that $L'$ is a $M_{LVE}$-dominating set of $G$. Hence, we have $\gamma(G)\leq \gamma(G')$. Thus, $\gamma(G)=\gamma(G')$.
\end{proof}

Observe that, from Lemma \ref{lem:5}, it follows that for every pair of edges $(e,e')$ where $e$ is an edge of an end block of $B'$ and $e'$ is an edge incident to the cut vertices $c_q$ and $c_i$ with $r(N_G(c_i)\setminus V(B'))=\emptyset$, we have  because $\max\{k(e)+k(e')-1,0\}\leq s'(c_i)+|r(N_G[c_q])|$. Therefore, condition $(iv)$ of $M_{LVE}$-domination holds for every such pair of edges. In the next lemma, we show that condition $(iv)$ holds for the pair of edges where one edge is from an end block of $B'$ and another is an edge incident to $c_h$ and $c_q$.

\begin{lem}\label{lem:6}
	Let $B'$ be a support block in $\sigma_{\mathcal{B}}$, $c_i\in C\setminus\{c_h\}$ is a cut vertex of $B'$ and $B_{ij}$ is an end block containing the vertex $c_i$. Also, assume that the labels of $G$ are as follows: $k(e)\leq |r(N_G[c_i])|$ for every edge $e$ of $B_{ij}$, $\max\{k(e)+k(e')-1,0\}\leq |r(N_G[c_i])|$ for every pair of edges $e,e'$ of $B_{ij}$, $s(v)\leq |r(N_G[c_i])|$ for every non-cut vertex $v$ of $B_{ij}$, and $\max\{k(e)+k(e')-1,0\}\leq |r(N_G[c_a]\cup N_G[c_b])|$ for every pair of edges $(e,e')$ such that $e$ and $e'$ are the edges of the end blocks $B_{aj}$ and $B_{bj'}$, where $B_{aj}$ and $B_{bj'}$ are end blocks containing the cut vertices $c_a$ and $c_b$, respectively. Let $q$ be the index of a cut vertex in $C\setminus\{c_h\}$ such that $|r(N_G[c_q])|$ is minimum of all $|r(N_G[c_i])|$ corresponding to every cut vertices $c_i\in C\setminus\{c_h\}$ having $l_2(c_i)=1$. Then, $\gamma(G)=\gamma(G')$, where $G'$ is obtained by relabelling $s(c_h)$ as $s'(c_h)=\max\{s(c_h),3-|r(N_G[c_q])|\}$ and every other label remains the same.
\end{lem}

\begin{proof}
	Let $L$ be a minimum $M_{LVE}$-dominating set of $G$. Since $l_2(c_q)=1$, there is only one edge, say $e_q=c_qv_x$, where $v_x$ is a vertex of an end block incident to $c_q$, with $k$-label $2$. Now, for every pair of edges having $k$-label $2$ and incident to $c_q$, $L$ contains at least $3$ vertices to satisfy the condition $(iv)$ of $M_{LVE}$-dominating set. Consider the edge $e'=c_hc_q$. Clearly, $k$-label of $e'$ is $2$ and hence $L$ contains at least $3$ vertices from $N_G[c_q]\cup N_G[c_h]$. Since $L$ contains every vertex with $t$-label $R$, $r(N_G[c_q]\cup N_G[c_h])$ is contained in $L$. Therefore, $L$ contains at least $3-(|r(N_G[c_q])|)$ vertices of $N_G[c_q]\cup N_G[c_h]$ having $t$-label $B$. Also, $L$ contains at least $s(v)$ vertices from $N_G[v]$ for every vertex $v\in V(G)$. Hence, there are $s(c_q)$ and $s(c_h)$ vertices in $L$ for the vertices $c_q$ and $c_h$, respectively. If $s(c_h)\geq 3-(|r(N_G[c_q])|)$, then $s'(c_h)=s(c_h)$ and thus $L$ is a $M_{LVE}$-dominating set of $G'$. Otherwise, suppose there are $3-(|r(N_G[c_q])|)$ vertices in $(N_G[c_q]\setminus V(B'))\cap L$ with $t$-label $B$. Let $L'$ be the set obtained from $L$ by replacing the $3-(|r(N_G[c_q])|)$ vertices of $N_G[c_q]\setminus V(B')$ having $t$-label $B$ with the $3-(|r(N_G[c_q]|)$ vertices of $N_G[c_h]$ having $t$-label $B$. Hence, $L'$ contains at least $\max\{s(c_h),3-(|r(N_G(c_q)\setminus V(B'))|)+|r(V(B'))|$ many vertices for the vertex $c$. Therefore, it is easy to observe that $L'$ is a $M_{LVE}$-dominating set of $G'$. Hence, we have $\gamma(G')\leq \gamma(G)$. 
	
	Conversely, let $L'$ be a minimum $M_{LVE}$-dominating set of $G'$. Now, every edge of $G'$ has the same $k$-label as $G$ and every vertex of $G'$ has the same $M_V$-label as $G$ except for the vertex $c_h$. We show that $L'$ is a $M_{LVE}$-dominating set of $G$. Since $L'$ contains at least $s'(c_h)=\max\{s(c_h),3-(|r(N_G(c_q)\setminus V(B'))|)+|r(V(B'))|\}$ vertices in $N_G[c_h]$, it is straightforward that $L'$ contains at least $s(c_h)$ vertices in $N_G[c_h]$. Since every other label remains the same, it is easy to see that $L'$ is a $M_{LVE}$-dominating set of $G$. Hence, we have $\gamma(G)\leq \gamma(G')$. Thus, $\gamma(G)=\gamma(G')$.
\end{proof}

From Lemma \ref{lem:6}, for the pair of edge $(e,e')$ where $e$ is an edge incident to $c_q$ with $k(e)=2$ and $e'=c_qc_h$, we have $\max\{k(e)+k(e')-1,0\}\leq s(c_h)+r(N_G[c_q])$. Therefore, the condition $(iv)$ of $M_{LVE}$-domination holds for every pair of edges incident to the cut vertices of $B'$.

\begin{lem}\label{lem:7}
	Let $B'$ be a support block in $\sigma_{\mathcal{B}}$, $c_i\in C\setminus\{c_h\}$ is a cut vertex of $B'$ and $B_{ij}$ is an end block containing the vertex $c_i$. Also, assume that the labels of $G$ are as follows: $k(e)\leq |r(N_G[c_i])|$ for every edge $e$ of $B_{ij}$, $\max\{k(e)+k(e')-1,0\}\leq |r(N_G[c_i])|$ for every pair of edges $e,e'$ of $B_{ij}$, $s(v)\leq |r(N_G[c_i])|$ for every non-cut vertex $v$ of $B_{ij}$, and $\max\{k(e)+k(e')-1,0\}\leq |r(N_G[c_a]\cup N_G[c_b])|$ for every pair of edges $(e,e')$ such that $e$ and $e'$ are the edges of the end blocks $B_{aj}$ and $B_{bj'}$, where $B_{aj}$ and $B_{bj'}$ are end blocks containing the cut vertices $c_a$ and $c_b$, respectively. Moreover, $\max\{k(e'')+k(e''')-1,0\}\leq r(N_G[c_i])+s(c')$ for every pair of edges $(e'',e''')$ such that $e''$ is an edge of an end block $B_{ij}$ and $e'''$ is an edge of $B'$ having endpoint in $c'$, where $B_{ij}$ contains the cut vertex $c_i$. Then, $\gamma(G)=\gamma(G')+|r(N_G[c_i]\setminus V(B'))|$, where $G'$ is obtained by deleting all the vertices of $N_G[c_i]\setminus V(B')$ and relabelling $k(e)$ and $s(c_i)$ as $k'(e)=\max\{k(e)-|r(N_G(c_i)\setminus V(B'))|,0\}$, for every edge $e$ of $B'$ incident to $c_i$ and $s'(c_i)=\max\{s(c_i)-|r(N_G(c_i)\setminus V(B'))|,0\}$, for the cut vertex $c_i$ of $B'$, respectively. Every other label remains the same.
\end{lem}

\begin{proof}
	Let $L$ be a minimum $M_{LVE}$-dominating set of $G$. Thus, $L$ contains at least $k(e)$ vertices of $N_G[e]$ for every edge $e$, and at least $s(v)$ vertices of $N_G[v]$ for every vertex $v$. Clearly, for the vertex $c_i$, $L$ contains at least $s(c_i)$ vertices of $N_G[c_i]$. Also, every vertex of $N_G[c_i]$ with $t$-label $R$ is in $L$. Therefore, $r(N_G[c_i])$ is contained in $L$. Let $L'$ be the set obtained from $L$ by replacing all the of $s(c_i)-|r(N_G[c_i])|$ vertices of $N_G[c_i]\setminus V(B')$ having $t$-label $B$ with all the of $s(c_i)$ vertices of $N_G[c_i]\setminus V(B')$ with $t$-label $B$. It is easy to verify that $L'$ is a minimum $M_{LVE}$-dominating set of $G$. Consider the set $L''=L'\setminus (r(N_G[c_i]\setminus V(B')))$. We show that $L'$ is a $M_{LVE}$-dominating set of $G'$. For every edge $e$ incident to $c_i$ in $G'$, $L''$ contains at least $\max\{k(e)-|r(N_G(c_i)\setminus V(B'))|,0\}=k'(e)$ vertices of $N_{G'}[e]$ because $L'$ contains at least $k(e)$ vertices of $N_G[e]$, where $N_{G'}[e]=N_G[e]\setminus (V(B')\setminus\{c_i\})$. Similarly, for the vertex $c_i$, $L''$ contains at least $\max\{s(c_i)-|r(N_G(c_i)\setminus V(B'))|,0\}=s'(c_i)$ vertices of $N_{G'}[c_i]$, where $N_{G'}[c_i]=N_G[c_i]\setminus (V(B')\setminus\{c_i\})$. Since every other label remains the same, $L''$ is a $M_{LVE}$-dominating set of $G'$. Hence, $\gamma(G')\leq \gamma(G)-|r(N_G(c_i)\setminus V(B'))|$.
	
	Conversely, Let $L'$ be a minimum $M_{LVE}$-dominating set of $G'$. Hence, $L'$ contains $s'(c_i)$ vertices of $N_{G'}[c_i]$ and $k'(e)$ vertices of $N_G[e]$ for every edge $e$ incident to $c_i$. Note that, all the other $M_V$-labels of the vertices in $V(G')$ and $k$-label of the edges are the same as $G$. Consider the set $L=L\cup (r(N_G[c_i]\setminus V(B')))$. We show that $L$ is a $M_{LVE }$-dominating set of $G$. Since $L'$ contains $s'(c)=\max\{s(c_i)-|r(N_G(c_i)\setminus V(B'))|,0\}$ vertices of $N_{G'}[c_i]$, we have $s'(c)+|r(N_G(c_i)\setminus V(B'))|=s(c_i)$ vertices of $N_{G}[c_i]$ in $L$. Also, $L'$ contains at least $k'(e)=\max\{k(e)-|r(N_G(c_i)\setminus V(B'))|,0\}$ vertices of $N_{G'}[e]$. This implies that $L$ contains $k'(e)+|r(N_G(c_i)\setminus V(B'))|=k(e)$ vertices of $N_G[e]$. Since every other label is the same and from the given conditions on the $k$-label and $s$-label of the edges and vertices in $G$, it is easy to verify that $L$ is a $M_{LVE}$-dominating set of $G$. Hence, we have $\gamma(G)\leq \gamma(G')+|r(N_G(c_i)\setminus V(B'))|$. Thus, we have $\gamma(G)=\gamma(G')+|r(N_G(c_i)\setminus V(B'))|$. 
\end{proof}

Next, we prove the correctness of the Algorithm \ref{Algo:StarBlock}. Note that, the lines $1-8$ and $9-11$ of Algorithm \ref{Algo:StarBlock} are similar to lines $6-18$ and $20-22$ of Algorithm \ref{Algo:Blockmain}, respectively. Hence, the correctness of these lines of Algorithm \ref{Algo:StarBlock} follows from Lemma \ref{lem:3if} and Lemma \ref{lem:2if}, respectively.

\begin{lem}\label{lem:8}
	Let $G$ be a star block graph having cut vertex $c'$. Also, assume that for every edge $e$ of an end block $B_{ij}$ which is not incident to $c'$, $k(e)\leq r(N_G[c'])$ and for any pair of edges $(e,e')$ of an end block $B_{ij}$ which is not incident to $c'$, $\max\{k(e)+k(e')-1,0\}\leq |r(N_G[c'])|$. Then $\gamma(G)=\gamma(G')$, where $G'$ is obtained by relabelling $\max\{s(c')-|r(N_G[c'])|,0\}$ vertices of $N_G[c']$ as $R$ and every other label remains the same.
\end{lem}

\begin{proof}
	Let $L$ be a minimum $M_{LVE}$-dominating set of $G$. Hence, $L$ contains at least $s(c')$ vertices of $N_G[c']$. Since every vertex having $t$-label $R$ is in $L$, we have $r(N_G[c'])$ is contained in $L$. Therefore, there are $\max\{s(c')-|r(N_G[c'])|,0\}$ many vertices of $N_G[c']$ having $t$-label $B$ is in $L$. Let $L'$ be a set obtained from $L$ by replacing the vertices of $N_G[c']$ having $t$-label $B$ with the vertices of $N_G[c']$ having $t'$-label $R$. Since every other label remains the same, it is easy to observe that $L'$ is a $M_{LVE}$-dominating set of $G'$. Hence, we have $\gamma(G')\leq \gamma(G)$.
	
	Conversely, let $L'$ be a minimum $M_{LVE}$-dominating set of $G'$. Since every edge of $G'$ has same $k$-label as $G$ and every vertex of $G'$ has the same $M_V$-label as $G$ except for $\max\{(s(c')-|r(N_G[c'])|),0\}$ vertices of $N_{G'}[c']=N_G[c']$. Therefore, $L'$ is a $M_{LVE}$-dominating set of $G$. Hence, we have $\gamma(G)\leq \gamma(G')$. Thus, $\gamma(G)=\gamma(G')$.
\end{proof}

\subsubsection{Running time analysis}
Now, we discuss the running time of our algorithm. The cut tree can be constructed in $O(n+m)$ time. In each iteration of end block, the algorithm takes $O(|E(B_{ij})|)$ time to find $l_2(B_{ij})$ and $l_1(B_{ij})$. Each condition of line $8,11,14$ can be checked in $O(1)$ time. $s_{max}$ can be determined by finding the maximum $s$-label of the vertices in $B_{ij}$ except for $c_i$. This can be done in $O(|V(B_{ij})\setminus \{c_i\}|)$ time. Hence, line $17$ can be executed in $O(|V(B_{ij})\setminus \{c_i\}|)$ time, that is $O(deg(c_i))$ time. Therefore, each iteration of for loop at line $6$ takes $O(|E(B_{ij})|)+O(deg(c_i))$ time. Thus, the total running time of the for loop at line $6-18$ takes $$O(\sum_{B_{ij}}|E(B_{ij})|)+O(deg(c_i))$$ time and the total running time of the for loop at line $5-18$ takes $$\sum_{c_{i}\in V(B')}(\sum_{B_{ij}}{O(|E(B_{ij})|)}+O(deg(c_i)))$$ time. Now, $l_2(c_i)$, $max_1(c_i)$, $max_2(c_i)$, and $r(N_G[c_i])$ can be determined in $O(deg(c_{i}))$ time. Also, the condition at line $21$ can be checked in $O(1)$ time. Therefore, the for loop at line $19-22$ takes $O(\sum_{c_i\in B'}deg(c_{i}))$ time. The conditions in line $23$ can be checked in $O(deg(c_h))$ time. Therefore, the lines $23-24$ can be executed in $O(deg(c_h))$ time. Finding the index $q$ takes $O(deg(c_h))$ time. Hence, line $25$ takes $O(deg(c_h))$ time. Now, the conditions of lines $26$ and $28$ can be checked in $\sum_{c_i\in V(B')}O(deg(c_i))$ time. Also, line $30$ takes $O(1)$ time. Therefore, the lines $5-30$ can be executed in $$\sum_{c_{i}\in V(B')}(\sum_{B_{ij}}{O(|E(B_{ij})|)}+O(deg(c_i)))$$ time. Inclusion of the required vertices in the already constructed set $L$ takes $O(\sum_{c_i\in V(B')}deg(c_i))$ time. Hence, line $31$ takes $\sum_{c_i\in V(B')}O(deg(c_i))$ time. Deletion of every end block vertices of $N_G(c_i)$ of a support block $B'$ takes $O(deg(c_i))$ time. Also, relabelling of the edges of $B'$ incident to $c_i$ can be done in $O(deg(c_i))$ time. Hence, lines $33-34$ can be executed in $O(\sum_{c_i\in V(B')}deg(c_i))$ time. Also, updating the $s$-label of every cut vertex of $B'$ can be done in $O(deg(c_h))$ time. Thus, lines $35-36$ take $O(deg(c_h))$ time. Therefore, the total running time of the for loop at line $4$ takes $$\sum_{B'\in \sigma_{\mathcal{B}}}O(\sum_{c_i\in V(B')}deg(c_i))+O(\sum_{c_i\in V(B')}\sum_{B_{ij}}{O(|E(B_{ij})|)}+O(deg(c_i)))+O(deg(c_h))$$ time. Observe that the StarBlock() algorithm has the same steps as the first two rounds of the main algorithm and thus, has the same running time as the first two rounds. Hence, the total running time of the algorithm is $$\sum_{B'\in \sigma_{\mathcal{B}}}O(\sum_{c_i\in V(B')}deg(c_i))+O(\sum_{c_i\in V(B')}\sum_{B_{ij}}{O(|E(B_{ij})|)+O(deg(c_i))})+O(deg(c_h))$$, that is, $O(|E(G)|)+O(|V(G)|)$. Therefore, a minimum liar's ve-dominating set can be computed in $O(n+m)$ time. 

\begin{theo}
	The \textsc{MinLVEDP} can be solved in linear time for block graphs.
\end{theo}

\subsection{Algorithm in proper interval graphs}

In this section, we describe a linear time algorithm to compute a minimum liar's ve-dominating set in connected proper interval graphs. Let $G=(V, E)$ be a connected proper interval graph. Also, assume that $\mathcal{I}=\{I(u):u\in V\}$ be the interval representation of $G$. Observe that, since $G$ is a proper interval graph, for any two distinct vertices $u,v$, we have $I(u)\nsubseteq I(v)$ and $I(v)\nsubseteq I(u)$.
\subsubsection{Description of the algorithm}
To describe the algorithm, we first define an ordering of the edges of $G$. Let us assume that every edge $e\in E(G)$ is represented as $e=uv$ such that $r(u)<r(v)$. For a pair of edges $e_i=u_iv_i$ and $e_j=u_jv_j$, we say $e_i<_E e_j$ if either $r(v_i)<r(v_j)$ or $v_i=v_j$ and $r(u_i)<r(u_j)$. An ordering $\sigma_E=(e_1, e_2, \ldots, e_m)$ of the edges of an interval graph $G=(V,E)$ is called an \emph{edge interval ordering} if $e_i<_E e_j$ for $i<j$. Let $S=\{v_1,v_2,\ldots,v_k\}$ be a subset of $V$ such that $r(v_k)<r(v_{k-1})< \ldots <r(v_1)$. For an integer $q$, let $T^q[S]$ denote the set of $q$ vertices from $S$ whose right endpoints are more. Therefore, $T^q[S]=\{v_1,v_2,\ldots,v_q\}$. 

\begin{algorithm}[h]
	\scriptsize
	\caption{LVEDS-PROPER-INTERVAL$(G)$}
	\textbf{Input:} A connected proper interval graph $G=(V,E)$.\\
	\textbf{Output:} A minimum liar's ve-dominating set $L$ of $G$.
	\begin{algorithmic}[1]
		\State Construct an interval representation of $G$;
		\State Find an edge interval ordering $\sigma_E=(v_1,v_2,\ldots,v_n\}$ of $E(G)$;
		\State $L=\emptyset$;
		\For{$(i=1~ to ~ m)$}
		\If{$(|N_G[e_i]\cap L|\leq 2)$}
		\State $q=2-|N_G[e_i]\cap L|$;
		\State $L=L\cup T^q[N_G[e_i]\setminus L]$;
		\State $min_1=\min\{k: k>i ~and ~ |N_G[e_k]\cap L|<2\}$;
		\State $min_2=\min\{k:k>i ~and~ N_G[e_k]|\cap L=N_G[e_i]\cap L\}$;
		\State $min=\min\{min_1,min_2\}$;
		\If{($min=min_1$)}
		\State $q_1=2-|N_G[e_{min}]\cap L|$;
		\State $L=L\cup T^{q_1}[N_G[e_{min}]\setminus L]$;
		\Else
		\State $L=L\cup T^{1}[N_G[e_{min}]\setminus L]$;
		\EndIf	
		\EndIf
		\State $i=min$;
		\EndFor
		\State \Return $L$;
	\end{algorithmic}
	\label{Algo:Interval}
\end{algorithm}

Now we give a brief description of our proposed algorithm. First, we order the edges according to edge interval ordering $\sigma_E$. In each iteration, we process an edge $e_i$ from $\sigma_E$. Also, we start with an empty set $L$. Let $E_i=\{e_1, e_2, \ldots, e_i\}$ denote the first $i$ edges of $\sigma_{E}$. The algorithm is an iterative algorithm. The idea is to select vertices in the set $L$ such that after $i$-th iteration, every edge in $E_i$ satisfies the conditions of liar's ve-domination. Each iteration consists of two rounds. At the first round (lines $3-7$), based on the cardinality of $N_G[e_i]\cap L$, we decide to include a few vertices in $L$ such that condition $(i)$ of liar's ve-domination holds for the edge $e_i$. Note that, for the pair of edges $(e_i,e_j)$, if $N_G[e_i]\cap N_G[e_j]= \emptyset$, then the condition $(ii)$ of liar's ve-domination follows from condition $(i)$ of liar's ve-domination. Hence, it is sufficient to check the condition $(ii)$ for the pair of edges $(e_i,e_j)$ such that $N_G[e_i]\cap N_G[e_j]\neq \emptyset$.  In the second round (lines $8-15$), we find the minimum index of an edge $e_{min}\in E\setminus E_i$ such that either $|N_G[e_{min}]\cap L|<2$ or $N_G[e_{min}]\cap L=N_G[e_i]\cap L$ and include some vertices of $N_G[e_{min}]$ in $L$. This step ensures that the condition $(ii)$ of liar's ve-domination holds for every pair of edges where one edge is $e_i$ and another is any edge $e_j\in E(G)\setminus\{e_i\}$ with $N_G[e_i]\cap N_G[e_j]\neq \emptyset$. Next, we skip all the edges between $e_i$ and $e_{min}$ in $\sigma_{E}$ and proceed to the next iteration. The detailed algorithm is described in Algorithm \ref{Algo:Interval}.

\subsubsection{Correctness}
In this subsection, we prove the correctness of our proposed algorithm LVEDS-PROPER-INTERVAL$(G)$. First, we introduce some notations that are used in proof of correctness. Now, the $i$-th iteration is executed in two rounds. Let $L_{i_1}$ and $L_{i}$ be the sets constructed by Algorithm \ref{Algo:Interval} after the first round and second round respectively. Initially, we set $L_0=\emptyset$. 
Also, for an edge $e_i= u_iv_i$, let $I(e_i)=[l(e_i),r(e_i)]=I(u_i)\cup I(v_i)$ denote the interval corresponding to the edge $e_i$. Thus, $r(e_i)=r(v_i)$ and $l(e_i)=l(u_i)$. The following observation is immediate from the interval representation of $G$. 

\begin{obs}\label{Obs:vnotintheedgee}
	At $i$-th iteration, let $e_i=u_iv_i$ be the edge that is going to be processed. If $v\notin N_G[e_i]$, then either $r(v)<l(e_i)$ or $l(v)>r(e_i)$.
\end{obs}


The following proposition guarantees that after the first round of processing of the edge $e_i$, the construction of the set $L_{i_1}$ is such that every edge $e_j$, where $i<j<min$, satisfies the conditions of liar's ve-domination.

\begin{pre}\label{Prop:everyedgefrmoitomin}
	At $i$-th iteration, after the first round of processing of the edge $e_i$, for any edge $e_j$ with $i<j<min$, $|N_G[e_j]\cap L_{i_1}|\geq 3$.
\end{pre}

\begin{proof}
	We prove this using contradiction. Let $e_j$ be an edge with $i<j<min$, such that $|N_G[e_j]\cap L_{i_1}|\leq 2$. Note that, if $|N_G[e_j]\cap L_{i_1}|<2$, then there is an edge $e_j$ with $i<j<min$ which contradicts the fact that $e_{min}$ be the minimum index edge with such property. Therefore, assume that $|N_G[e_j]\cap L_{i_1}|=2$. Since $j<min$, we have $N_G[e_j]\cap L_{i_1}\neq N_G[e_i]\cap L_{i_1}$. Thus, there is a vertex $v\in N_G[e_j]\cap L_{i_1}$ such that $v\notin N_G[e_i]$. Hence, by Observation \ref{Obs:vnotintheedgee} we have $r(e_i)<l(v)$ or $l(e_i)>r(v)$. Since $v\in N_G[e_j]$, we have $r(e_j)>l(v)$ and $l(e_j)<r(v)$. Thus, we have $l(e_j)<r(v)<l(e_i)$ and using the fact that $r(e_j)\geq r(e_i)$ as $j>i$, we have $I(e_i)\subseteq I(e_j)$. Thus, for any vertex $u\in N_G[e_i]$, we have $I(u)\cap I(e_j)\neq \emptyset$. This implies that $|N_G[e_j]\cap L_{i_1}|\geq 3$ which is a contradiction. Therefore, after the first round, for any edge $e_j$ with $i<j<min$, we have $|N_G[e_j]\cap L_{i_1}|\geq 3$. 
\end{proof}

The next two propositions ensure that after the second round of processing of the edge $e_i$, every pair of edges $(e_i,e_j)$, where $e_j$ is an edge in $\sigma_E$ with $j>min$, and $N_G[e_i]\cap N_G[e_j]\neq \emptyset$, satisfies the condition $(ii)$ of liar's ve-domination.

\begin{pre}\label{Prop:secondcondition}
	At $i$-th iteration, after the first round of processing of the edge $e_i$, if $|N_G[e_{min}]\cap L_{i_1}|<2$, then for every edge $e_j$ such that $j>min$ and $N_G[e_i]\cap N_G[e_j]\neq \emptyset$ , there is a vertex $v\in T^q[N_G[e_{min}]\setminus L_{i_1}]$ such that $v\in N_G[e_i]\cup N_G[e_j]$, where $q=2-|N_G[e_{min}]\setminus L_{i_1}|$. 
\end{pre}

\begin{proof}
	Let $v\in T^q[N_G[e_{min}]\setminus L_{i_1}]$ be a vertex such that $v\notin N_G[e_i]$. Hence, by Observation \ref{Obs:vnotintheedgee} either $l(v)>r(e_i)$ or $r(v)<l(e_i)$. Also, note that $v\in N_G[e_{min}]$ and $min>i$. Thus, $r(v)>l(e_{min})$ and $r(e_{min})>l(v)$. Suppose that $r(v)<l(e_i)$. Therefore, we have $I(e_i)\subset I(e_{min})$ which contradicts the fact that $|N_G[e_{min}]\cap L_{i_1}|<2$. Hence,  $l(v)>r(e_i)$. Since $N_G[e_i]\cap N_G[e_j]\neq \emptyset$, suppose that $u\in N_G[e_i]\cap N_G[e_j]$. Now, $l(u)<r(e_j)$ and $r(u)>l(e_j)$. Also, $l(u)<r(e_i)$ and $r(u)>l(e_i)$. Therefore, $l(u)<r(e_i)<l(v)$ and if $r(u)>r(v)$, then $I(v)\subset I(u)$ which is a contradiction. Hence, $r(u)<r(v)$. Now, $r(e_j)>r(e_{min})>l(v)$ and $r(v)>r(u)>l(e_j)$. Therefore, $I(v)\cap I(e_j)\neq \emptyset$. Hence, $v\in N_G[e_i]\cup N_G[e_j]$.
\end{proof}

\begin{pre}\label{Prop:thirdcondition}
	At $i$-th iteration, after the first round of processing of the edge $e_i$, if $N_G[e_{min}]\cap L_{i_1}= N_G[e_i]\cap L_{i_1}$, then for every edge $e_j$ such that $j>min$ and $N_G[e_i]\cap N_G[e_j]\neq \emptyset$, there is a vertex $v\in T^1[N_G[e_{min}]\setminus L_{i_1}]$ such that $v\in N_G[e_i]\cup N_G[e_j]$.
\end{pre}

\begin{proof}
	The proof is similar to the proof of Proposition \ref{Prop:secondcondition} and hence, omitted.
\end{proof}
Note that, after the second round of processing of $e_i$, either $L_{i}=L_{i_1}\cup T^q[N_G[e_{min}]\setminus L_{i_1}]$ or $L_{i}=L_{i_1}\cup T^1[N_G[e_{min}]\setminus L_{i_1}]$. Therefore, we have $|(N_G[e_i]\cup N_G[e_j])\cap L_i|\geq 3$, for any edge $e_j$ with $j\geq min$, and $N_G[e_j]\cap N_G[e_i]\neq \emptyset$.

In Algorithm \ref{Algo:Interval}, at $i$-th iteration there are two cases. If $|N_G[e_i]\cap L_{i-1}|\geq 3$, then $L_{i-1}=L_i$. If $|N_G[e_i]\cap L_{i-1}|\leq 2$, then it is easy to see that after the first round of $i$-th iteration, we have $|N_G[e_i]\cap L_{i_1}|= 2$. Since $L_{i_1}\subseteq L_i$, for every edge $e_j\in E_{i}$, we have $|N_G[e_i]\cap L_{i}|\geq 2$. Hence, we have the following observation.

\begin{obs}\label{Obs:firstcondition}
	At $i$-th iteration after the first round of processing of $e_i$, $L_{i}$ satisfy the condition that for every $e_j\in E_{i}$, $|N_G[e_j]\cap L_{i}|\geq 2$.
\end{obs}

Now, we show that after $i$-th iteration, every edge in $E_i$ satisfies the conditions of liar's ve-domination. 

\begin{pre}\label{Prop:Liar'sdomset}
	After the first round of processing of $e_i$, let $L_{i_1}$ be the set such that $|N_G[e_i]\cap L_{i_1}|\geq 2$. Also, assume that $e_{min}$ be the edge with minimum index in $\sigma_E$ such that $min>i$ and either $|N_G[e_{min}]\cap L_{i_1}|<2$ or $N_G[e_{min}]\cap L_{i_1}=N_G[e_i]\cap L_{i_1}$, where $|N_G[e_i]\cap L_{i_1}|=2$. After $i$-th iteration, $L_{i}$ satisfies the following properties:
	
	\noindent(a) $|N_G[e_k]\cap L_{i}|\geq 2$, for every edge $e_k$ in $E_i$
	
	\noindent(b) $|(N_G[e_k]\cup N_G[e_a])\cap L_{i}|\geq 3$, for every pair of distinct edges $e_k,e_a\in E_i$
	
	\noindent(c) $|(N_G[e_k]\cup N_G[e_a])\cap L_{i}|\geq 3$, where $e_k\in E_i$ and $e_a\in E\setminus E_i$ such that either $i<a\leq min$ or $a>min$ and $N_G[e_i]\cap N_G[e_a]\neq \emptyset$.
\end{pre}

\begin{proof}
	We prove this using induction on $i$. The proof of part $(a)$ follows from Observation \ref{Obs:firstcondition}. 
	
	Next, we prove part $(b)$. At $i=1$, there is only one edge in $E_i$ and thus, the statement is vacuously true. After the second round of processing of the edge $e_{i-1}$, let $L_{(i-1)_1}$ be the set of vertices such that $|N_G[e_{i-1}]\cap L_{(i-1)_1}|\geq 2$ and $e_{{min}_1}$ be the edge with minimum index in $\sigma_E$ such that ${min}_1>(i-1)$ and either $|N_G[e_{{min}_1}]\cap L_{(i-1)_1}|\geq 2$ or $N_G[e_{{min}_1}]\cap L_{(i-1)_1}=N_G[e_{i-1}]\cap L_{(i-1)_1}$. After $(i-1)$-th iteration, let $L_{i-1}$ be a set having the following properties: for every edge $e_k\in E_{i-1}$, $|N_G[e_k]\cap L_{i-1}|\geq 2$, for every pair of distinct edges $e_k,e_a\in E_{i-1}$, $|(N_G[e_k]\cup N_G[e_a])\cap L_{i-1}|\geq 3$, and for every pair of edges $(e_k,e_a)$, where $e_k\in E_{i-1}$ and $e_a$ be an edge such that either $(i-1)<a\leq min_1$ or $a>min_1$ and $N_G[e_{i-1}]\cap N_G[e_a]\neq \emptyset$. Observe that, $i\leq min_1$ and $L_{i-1}\subseteq L_i$. Therefore, for every edge $e_k\in E_{i-1}$, we have $|(N_G[e_k]\cup N_G[e_i])\cap L_i|\geq 3$.

	Finally, we prove part $(c)$. At $i=1$, after the first round $|N_G[e_1]\cap L_{1_1}|=2$, where $L_{1_1}=T^2[N_G[e_1]]$. Let $u\in T^2[N_G[e_1]]$ such that $r(u)\leq r(v)$ for any $v\in T^2[N_G[e_1]]$. Since $G$ is connected and $r(e_2)\leq r(u)$, we have $r(e_2)\leq r(v)$ and thus $l(e_2)\leq r(v)$ for any $v\in T^2[N_G[e_1]]$. Also, $l(v)<r(e_1)\leq r(e_2)$ and therefore, $I(v)\cap I(e_2)\neq \emptyset$ for any $v\in T^2[N_G[e_1]]$. Hence, $N_G[e_2]\cap L_{1_1}=N_G[e_1]\cap L_{1_1}$. Since $L_{1_1}\subseteq L_1$, using Proposition \ref{Prop:thirdcondition}, we have $|(N_G[e_1]\cup N_G[e_a])\cap L_{1}|\geq 3$ for any edge $e_a$ with $a\geq 2$, and $N_G[e_1]\cap N_G[e_a]\neq \emptyset$. Let us assume that after $(i-1)$-th iteration, $|(N_G[e_k]\cup N_G[e_a])\cap L_{i-1}|\geq 3$, for every pair of edges $e_k\in E_{i-1}$ and $e_a\in E\setminus\{e_k\}$ with $N_G[e_a]\cap N_G[e_k]\neq \emptyset$. From Observation \ref{Obs:firstcondition}, we have $|N_G[e_i]\cap L_{i_1}|\geq 2$. If $|N_G[e_i]\cap L_{i_1}|\geq 3$, then we are done. Otherwise, let $|N_G[e_i]\cap L_{i_1}|=2$ and $e_{min}$ be the minimum index edge such that $min>i$ and $N_G[e_i]\cap L_{i_1}=N_G[e_{min}]\cap L_{i_1}$ where $L_{i_1}=L_{i-1}\cup T^q[N_G[e_i]\setminus L_{i-1}]$. From Proposition \ref{Prop:secondcondition} and \ref{Prop:thirdcondition}, for every pair of edges $e_j$ with $j\geq min$ and $N_G[e_i]\cap N_G[e_j]\neq \emptyset$, we have $|(N_G[e_i]\cup N_G[e_{j}])\cap L_{i}|\geq 3$. Let $e_k$ be an edge such that $i<k<min$. From Proposition \ref{Prop:everyedgefrmoitomin}, we have $|N_G[e_k]\cap L_{i}|\geq 3$. Hence, $|(N_G[e_k]\cup N_G[e_a])\cap L_{i}|\geq 3$, for every pair of edges $e_k\in E_i$ and $e_a\in E\setminus E_i$ such that either $i<a\leq min$ or $a>min$ and $N_G[e_a]\cap N_G[e_i]\neq \emptyset$.
\end{proof}

Next, we show the minimality of the liar's ve-dominating set returned by Algorithm \ref{Algo:Interval}. We prove this using induction. We show that after $i$-th iteration, under different cases, $L_i$ is contained in a minimum liar's ve-dominating set. Therefore, at the end of the algorithm, it returns a minimum liar's ve-dominating set.

\begin{lem}
	At $i$-th iteration, let $|N_G[e_i]\cap L_{i-1}|\leq 2$. If $L_{i-1}$ is contained in some minimum liar's ve-dominating set $L$ of $G$, then there exists a minimum liar's ve-dominating set $L'$ of $G$ containing $L_i=L_{i-1}\cup T^q[N_G[e_i]\setminus L_{i-1}]$, where $q=2-|N_G[e_i]\cap L_{i-1}|$. 
\end{lem}

\begin{proof}
	Since $|N_G[e_i]\cap L_{i-1}|\leq 2$, there are $q=2-|N_G[e_i]\cap L_{i-1}|$ number of vertices, namely $\{v_1,v_2,\ldots,v_q\}$ in $N_G[e_i]\cap (L\setminus L_{i-1})$. Let us assume that $\{w_1,w_2\ldots,w_q\}\in T^q[N_G[e_i]\setminus L_{i-1}]$. If $v_j\in T^q[N_G[e_i]\setminus L_{i-1}]$ for every $j\in \{1,2,\ldots,q\}$, then we are done. Therefore, suppose that there is a vertex $v_j\in \{v_1,v_2,\ldots,v_q\}$ such that $v_j\notin T^q[N_G[e_i]\setminus L_{i-1}]$. Clearly, $r(v_j)<r(w_k)$, for any $k\in \{1,2,\ldots,q\}$. Let $e_a$ be an edge such that $v_j\in N_G[e_a]$ and $a\neq i$. Therefore, $r(v_j)>l(e_a)$ and $l(v_j)<r(e_a)$. Suppose that $a<i$. Now, from Proposition \ref{Prop:Liar'sdomset}, after $(i-1)$-th iteration, for any edge in $e_a\in E_{i-1}$, we have $|N_G[e_a]\cap L_{i-1}|\geq 2$ and for any edge $e_b$ such that $b>(i-1)$ and $N_G[e_a]\cap N_G[e_b]\neq \emptyset$, we have $|(N_G[e_a]\cup N_G[e_b])\cap L_{i-1}|\geq 3$. This implies that the vertices of $L_{i-1}$ is sufficient to satisfy the conditions of liar's ve-domination for the edge $e_a$. Also, if $a>i$, then $r(e_a)\geq r(e_i)$. Since $r(v_j)<r(w_k)$ for any $k\in\{1,2,\ldots,q\}$, we have $r(w_k)>l(e_a)$. Note that, $l(v_j)<l(w_k)$ otherwise $I(v_j)\subseteq I(w_k)$. Since $w_k\in N_G[e_i]$, we have $l(w_k)<r(e_i)\leq r(e_a)$. Hence $I(w_k)\cap I(e_a)\neq \emptyset$. Therefore, for every edge $e_a$ with $a>i$ which is ve-dominated by $v_j$, is also ve-dominated by $w_k\in T^q[N_G[e_i]\setminus L_{i-1}]$. Hence, $L'=(L\setminus\{v_1,v_2,\ldots,v_q\})\cup \{w_1,w_2,\ldots, w_q\}$ is a minimum liar's ve-dominating set of $G$. 
\end{proof}

\begin{lem}
	At $i$-th iteration, after the first round of processing of the edge $e_i$, let $L_{i-1}\subseteq L_{i_1}$, $|N_G[e_i]\cap L_{i_1}|=2$, and $e_{min}$ be the edge in $\sigma_E$ such that $min>i$, $|N_G[e_{min}]\cap L_{i_1}|<2$ and for every edge $i<p<min$, $|N_G[e_p]\cap L_{i_1}|\geq 3$. If $L_{i_1}$ is contained in some minimum liar's vertex-edge dominating set $L$ of $G$, then there exists a minimum liar's vertex-edge dominating set $L'$ of $G$ containing $L_i=L_{i_1}\cup T^q[N_G[e_{min}]\setminus L_{i_1}]$, where $q=2-|N_G[e_{min}]\setminus L_{i_1}|$.
\end{lem}

\begin{proof}
	Since $|N_G[e_{min}]\cap L_{i_1}|<2$, there are $q=2-|N_G[e_{min}]\setminus L_{i_1}|$ vertices, namely $\{v_1,v_2,\ldots,v_q\}$ in $N_G[e_{min}]\cap (L\setminus L_{i_1})$. Let us assume that $\{w_1,w_2\ldots,w_q\}\in T^q[N_G[e_{min}]\setminus L_{i_1}]$.  If $v_j\in T^q[N_G[e_i]\setminus L_{i-1}]$ for every $j\in \{1,2,\ldots,q\}$, then we are done. Therefore, suppose that there is a vertex $v_j\in \{v_1,v_2,\ldots,v_q\}$ such that $v_j\notin T^q[N_G[e_i]\setminus L_{i-1}]$. Clearly, $r(v_j)<r(w_k)$, for any $k\in \{1,2,\ldots,q\}$. Let $e_a$ be an edge such that $v_j\in N_G[e_a]$ and $a\neq i$. Therefore, $r(v_j)>l(e_a)$ and $l(v_j)<r(e_a)$. Suppose that $a<i$. Now, from Proposition \ref{Prop:Liar'sdomset}, after $(i-1)$-th iteration, for any edge in $e_a\in E_{i-1}$, we have $|N_G[e_a]\cap L_{i-1}|\geq 2$ and for any edge $e_b$ such that $b>(i-1)$ and $N_G[e_a]\cap N_G[e_b]\neq \emptyset$, we have $|(N_G[e_a]\cup N_G[e_b])\cap L_{i-1}|\geq 3$. This implies that the vertices of $L_{i-1}$, and hence vertices of $L_{i_1}$, are sufficient to satisfy the conditions of liar's ve-domination for the edge $e_a$. Now, for any edge $e_p$ with $i<p<min$, we have $|N_G[e_p]\cap L_{i_1}|\geq 3$. Therefore, if $i<a<min$, then $|N_G[e_a]\cap L_{i_1}|\geq 3$. Thus, vertices of $L_{i_1}$ are sufficient to satisfy the conditions of liar's ve-domination for the edge $e_a$. 
	Finally, assume that $a\geq min$. Since $r(v_j)<r(w_k)$ for any $k\in\{1,2,\ldots,q\}$, we have $r(w_k)>l(e_a)$. Note that, $l(v_j)<l(w_k)$, otherwise $I(v_j)\subseteq I(w_k)$. Since $w_k\in N_G[e_{min}]$, we have $l(w_k)<r(e_{min})\leq r(e_a)$. Hence $I(w_k)\cap I(e_a)\neq \emptyset$. Therefore, for every edge $e_a$ with $a>min$, which is ve-dominated by $v_j$, is also ve-dominated by $w_k\in T^q[N_G[e_{min}]\setminus L_{i_1}]$. Hence, $L'=(L\setminus\{v_1,v_2,\ldots,v_q\})\cup \{w_1,w_2,\ldots, w_q\}$ is a minimum liar's ve-dominating set of $G$.
\end{proof}

\begin{lem}
	At $i$-th iteration, after the first round of processing of the edge $e_i$, let $L_{i-1}\subseteq L_{i_1}$, $|N_G[e_i]\cap L_{i_1}|=2$, and $e_{min}$ be the edge in $\sigma_E$ such that $min>i$,  $N_G[e_{min}]\cap L_{i_1}=N_G[e_i]\cap L_{i_1}$, and for every edge $i<p<min$, $|N_G[e_p]\cap L_{i_1}|\geq 3$. If $L_{i_1}$ is contained in some minimum liar's vertex-edge dominating set $L$ of $G$, then there exists a minimum liar's vertex-edge dominating set $L'$ of $G$ containing $L_i=L_{i_1}\cup T^1[N_G[e_{min}]\setminus L_{i_1}]$.
\end{lem}

\begin{proof}
	Since $N_G[e_{min}]\cap L_{i_1}=N_G[e_i]\cap L_{i_1}$, to satisfy the condition $(ii)$ of liar's ve-domination for the pair of edges $(e_i,e_{min})$, there is a vertex, say $v$, in $(N_G[e_i]\cup N_G[e_{min}])\cap (L\setminus L_{i_1})$. Let us assume that $\{w\}= T^1[N_G[e_{min}]\setminus L_{i_1}]$. If $v=w$, then we are done. So, let $v\neq w$. Also, let $e_a$ is an edge such that $v\in N_G[e_a]$ and $a\neq i$. Therefore, $r(v)>l(e_a)$ and $l(v)<r(e_a)$. Suppose that $a<i$. Now, from Proposition \ref{Prop:Liar'sdomset}, after $(i-1)$-th iteration, for any edge in $e_a\in E_{i-1}$, we have $|N_G[e_a]\cap L_{i-1}|\geq 2$ and for any edge $e_b$ such that $b>(i-1)$ and $N_G[e_a]\cap N_G[e_b]\neq \emptyset$, we have $|(N_G[e_a]\cup N_G[e_b])\cap L_{i-1}|\geq 3$. This implies that the vertices of $L_{i-1}$, and hence vertices of $L_{i_1}$, are sufficient to satisfy the conditions of liar's ve-domination for the edge $e_a$. Now, for any edge $e_p$ with $i<p<min$, we have $|N_G[e_p]\cap L_{i_1}|\geq 3$. Therefore, if $i<a<min$, then $|N_G[e_a]\cap L_{i_1}|\geq 3$. Thus, vertices of $L_{i_1}$ are sufficient to satisfy the conditions of liar's ve-domination for the edge $e_a$. Finally, assume that $a\geq min$. Since $min>i$, $r(v)<r(w)$. Since $r(v)<r(w)$, we have $r(w)>l(e_a)$. Also, $l(v)<l(w)$ otherwise $I(v)\subseteq I(w)$. Since $w\in N_G[e_{min}]$, we have $l(w)<r(e_{min})\leq r(e_a)$. Hence $I(w_k)\cap I(e_a)\neq \emptyset$. Therefore, for every edge $e_a$ with $a>min$, which is ve-dominated by $v$, is also ve-dominated by $w$. Hence, $L'=(L\setminus\{v\})\cup \{w\}$ is a minimum liar's ve-dominating set of $G$.
\end{proof}

\subsubsection{Running time analysis}
Now, we discuss the running time of Algorithm \ref{Algo:Interval}.
%
Observe that, the vertices of $N_G[e_i]$ can be arranged according to the decreasing order of the right endpoints of intervals. Therefore, the vertices of $L_i$ are also arranged according to the decreasing order of the right endpoint of intervals. We can maintain the arranged order in $L_i$ while selecting a new vertex from $N_G[e_i]$ in $L_i$ by adding it to the beginning of $L_{i-1}$. More precisely, suppose that $L_{i-1}=\{v_1,v_2,\ldots,v_p\}$, where $r(v_1)>r(v_2)>\ldots>r(v_p)$, and a vertex $v_j$ is selected by the algorithm to be included in $L_i$. The following observation gives the condition of insertion in $L_i$.

\begin{obs}\label{Obs:InsertioninL}
	At the $i$-th iteration, let $T^q[N_G[e_i]\setminus L_{i-1}]$ be the set of vertices selected by Algorithm \ref{Algo:Interval} and $L_{i-1}=\{v_1,v_2,\ldots,v_p\}$. Then either $r(v_1)<r(v_j)$ or $r(v_2)<r(v_j)$ or $r(v_3)<r(v_j)$, where $v_j\in T^q[N_G[e_i]\setminus L_{i-1}]$. 
\end{obs}

\begin{proof}
	Suppose that $r(v_1)>r(v_j)$, $r(v_2)>r(v_j)$ and $r(v_3)>r(v_j)$ for some $v_j\in T^q[N_G[e_i]\setminus L_{i-1}]$. Since $v_j\in N_G[e_i]$, from interval representation of an edge, we have $I(v_j)\cap I(e_i)\neq \emptyset$. This implies that $r(v_j)>l(e_i)$ and $l(v_j)<r(e_i)$. Let us assume that there is a vertex $v\in L_{i-1}$ such that $r(v)>r(v_j)$. Hence, there is an edge $e_a\in E_{i-1}$ which is ve-dominated by $v$, that is $v\in N_G[e_a]$. Thus, we have $I(v)\cap I(e_a)\neq \emptyset$. Therefore, $r(v)>l(e_a)$ and $l(v)<r(e_a)$. Since $a<i$, we have $r(e_a)\leq r(e_i)$. As $r(v)>l(e_i)$ and $l(v)<r(e_a)\leq r(e_i)$, we have $I(v)\cap I(e_i)\neq \emptyset$. Therefore, $v\in N_G[e_i]$. Now, if $r(v)>r(v_j)$ for all $v\in \{v_1,v_2,v_3\}$, then we have $\{v_1,v_2,v_3\}\subseteq N_G[e_i]$. This implies that $|N_G[e_i]\cap L_{i-1}|\geq 3$ and thus $T^q[N_G[e_i]\setminus L_{i-1}]=\emptyset$. Thus, for any $v_j\in T^q[N_G[e_i]\setminus L_{i-1}]$, either $r(v_1)<r(v_j)$ or $r(v_2)<r(v_j)$ or $r(v_3)<r(v_j)$.  
	%
\end{proof}

\begin{obs}\label{Obs:checkingintersection}
	At $i$-th iteration, we can determine $|N_G[e_i]\cap L_{i-1}|\leq 3$ in $O(1)$ time.
\end{obs}

\begin{proof}
	By Observation \ref{Obs:vnotintheedgee}, we can determine whether a vertex of $L_{i-1}$ is in $N_G[e_i]$ or not in $O(1)$ time. Also, if $v\in L_{i-1}$ such that $v\notin N_G[e_i]$, then either $r(v)<l(e_i)$ or $l(v)>r(e_i)$. If $l(v)>r(e_i)$, then $l(v)>r(e_j)$ for any edge $e_j$ with $j<i$ and thus $v\notin N_G[e_j]$. This contradicts the fact that $v\in L_{i-1}$. Therefore, $r(v)<l(e_i)$. Since vertices of $L$ are in decreasing order of right endpoints of their interval representation, if $r(v)<l(e_i)$, then for any vertex $w\in L$ such that $r(w)<r(v)$, we have $r(w)<l(e_i)$ which implies that $w\notin N_G[e_i]$. Hence, to check whether $|N_G[e_i]\cap L_{i-1}|\leq 3$, it requires to check at most three vertices of $L_{i-1}$ is in $N_G[e_i]$ or not. Therefore, $|N_G[e_i]\cap L_{i-1}|\leq 3$ can be computed in $O(1)$ time.
\end{proof}
Given a proper interval graph $G=(V,E)$ with $|V|=n$ and $|E|=m$, we can construct the interval representation of $G$ in $O(n+m)$ time. Thus, we can sort the edge set $E$ in increasing order in $O(n+m)$ time. Hence, line $1$ and $2$ of Algorithm \ref{Algo:Interval} takes $O(n+m)$ each. Also, we can arrange the vertices of $N_G[v]$ for every $v\in V$ in descending order of $r(v)$ in $O(n+m)$ time using the interval representation. Now suppose that at $i$-th iteration, the edge $e_i=u_iv_i$ is being processed. Using Observation \ref{Obs:checkingintersection} we can determine $|N_G[e_i]\cap L_{i-1}|\leq 3$, in $O(1)$ time. Therefore, line $5$ takes $O(1)$ time. Based on the value of $q$, we add $q$ vertices in $L_{i-1}$. Note that, to construct the set $T^q[N_G[e_i]]$, we don't need to check the whole neighbourhood of $e_i=u_iv_i$. Observe that, we need only three vertices of $N_G[e_i]$ sorted according to the decreasing order of their right endpoint in interval representation. Now, vertices of $N_G[e_i]$ are the vertices of $N_G[u_i]$ and $N_G[v_i]$. Since vertices of $N_G[u_i]$ and $N_G[v_i]$ are already arranged in descending order, we can find the set $T^q[N_G[e_{i}]]$ in $O(1)$ time where $q\leq 2$. Now, at a time at most two vertices can be added to $L_{i-1}$. So, by Observation \ref{Obs:InsertioninL} we can maintain the ordering of $L_i$ in $O(1)$ time. Hence, line $7$ can be executed in $O(1)$ time. Let $L_{i_1}$ be the set obtained from $L_{i-1}$ by adding $T^q[N_G[e_i]\setminus L_{i-1}]$ vertices. Now, we determine $min$ by checking the $N_G[e_k]\cap L_{i_1}$, for the edges $e_k$ with $k>i$. As discussed above, this condition can be checked in $O(1)$ time for each edge $e_k$ with $i<k\leq min$. Therefore, $min$ can be determined in $O(min-i)$ time. Hence line $8-10$ takes $O(min-i)$ time. Also, the condition at line $11$ takes $O(1)$ time. We have already discussed that adding a vertex to $L$ takes $O(1)$ time. Therefore, line $13-15$ takes $O(1)$ time. Hence, the total running time of an iteration of the for loop at line $4$ is $O(min-i)$. Note that, the next iteration of the algorithm is updated to $min$. Thus, the total running time of the for loop at line $4$ is $O(m)$. Therefore, the running time of the algorithm is $O(n+m)$. Hence, we have the following theorem.

\begin{theo}
	\textsc{MinLVED} can be solved in linear time for proper interval graphs.
\end{theo}

\section{NP-completeness in undirected path graph}\label{sec:undirectedpath}
We show that the decision version of liar's ve-domination problem is NP-complete for \emph{undirected path graphs}. 
We prove the NP-completeness by reducing an instance of the $3$-dimensional matching problem, which is known to be NP-complete \cite{garey1979}, to an instance of \textsc{DecideLVEDP} for undirected path graphs. The $3$-dimensional matching problem is defined as follows:

\noindent\underline{\textbf{$3$-Dimensional Matching Problem} (\textsc{$3$-DMP})}

\noindent\emph{Instance}: A set $M\subseteq U\times V\times W$, where $U$, $V$ and $W$ are disjoint sets with $|U|=|V|=|W|=q$.

\noindent\emph{Question}: Does there exist a subset $M'$ of $M$ of size $q$ such that no two elements of $M'$ agree in any coordinates?

\begin{theo}
	The \textsc{DecideLVEDP} is NP-complete for the undirected path graph.
\end{theo}

\begin{proof}
	It is easy to see that \textsc{DecideLVEDP} is in NP. First, we construct a clique tree $\mathcal{T}$ of an undirected path graph $G$ from an instance of \textsc{$3$-DMP}. Let $U=\{u_i:1\leq i\leq q\}$, $V=\{V_i:1\leq i\leq q\}$, $W=\{w_i:1\leq i\leq q\}$ and $M=\{m_i=(u_j,v_k,w_s):u_j\in W,v_k\in V, w_s\in W~and~ 1\leq i\leq p\}$ be an instance of \textsc{$3$-DMP}. The vertices of $\mathcal{T}$ are as follows:

	\begin{figure}[h!]
		\centering
		\includegraphics[scale=0.5]{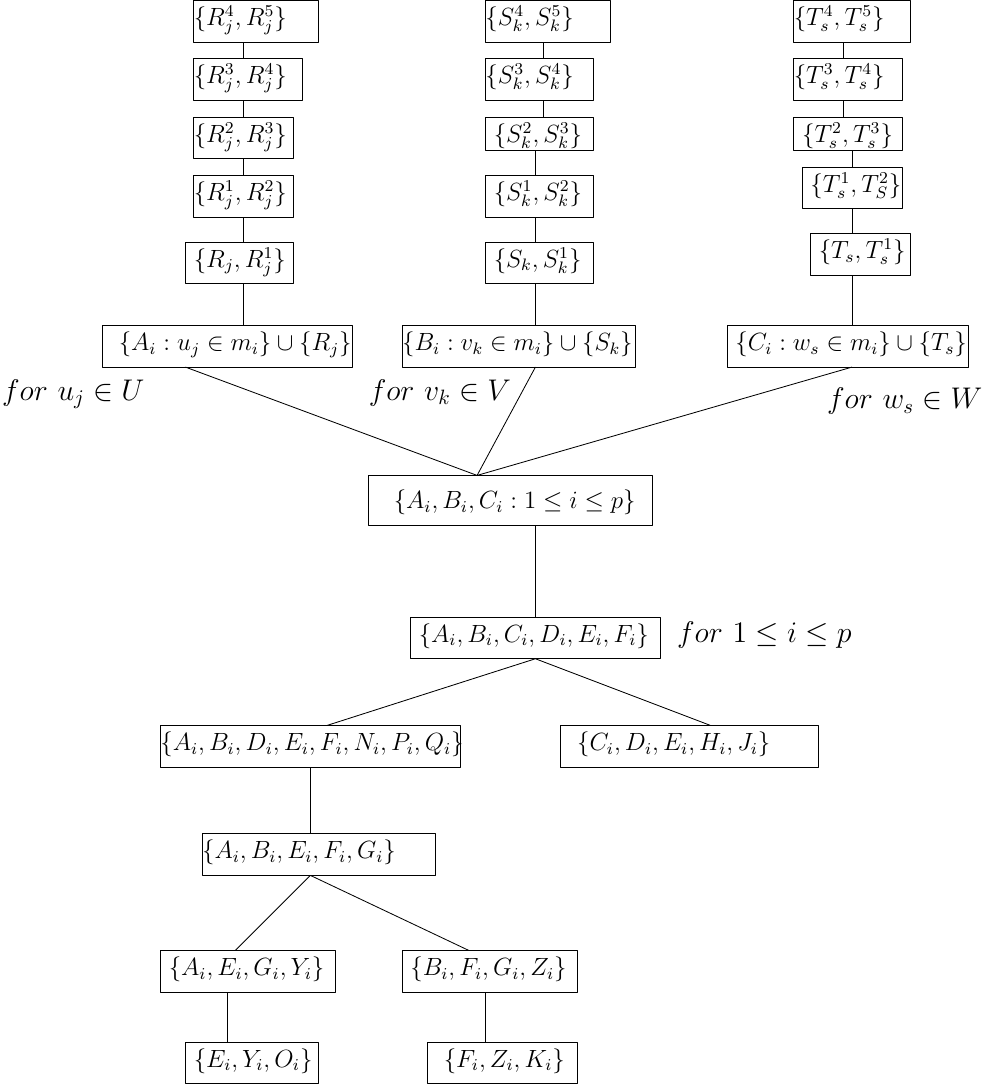}
		\caption{Clique tree $T$}
		\label{fig:Undirectedpath}
	\end{figure}
	
	For each $m_i\in M$, we consider a set of $8$ cliques, namely,  $\{A_i,B_i,C_i,D_i,E_i,\\F_i\}$, 
	$\{A_i,B_i,D_i,E_i,F_i,N_i,P_i,Q_i\}$, 
	$\{C_i,D_i,E_i,H_i,J_i\}$, 
	$\{A_i,B_i,E_i,F_i,G_i\}$, 
	$\{A_i\\,E_i,G_i,Y_i\}$, 
	$\{B_i,F_i,G_i,Z_i\}$, 
	$\{E_i,Y_i,O_i\}$, 
	$\{F_i,Z_i,K_i\}$. These cliques are corresponding to the triple $m_i$ and do not depend on the components of the triple. Now, for each $u_j\in U$, we consider a clique $(\{A_i:u_j\in m_i\}\cup \{R_j\})$ and also, five other cliques $\{R_j,R^1_j\}$, $\{R^1_j,R^2_j\}$, $\{R^2_j,R^3_j\}$, $\{R^3_j,R^4_j\}$ and $\{R^4_j,R^5_j\}$. Similarly, corresponding to each $v_k\in V$ and $w_s\in W$, we consider the six cliques $(\{B_i:v_k\in m_i\}\cup \{S_k\})$, $\{S_k,S^1_k\}$, $\{S^1_k,S^2_k\}$, $\{S^2_k,S^3_k\}$, $\{S^3_k,S^4_k\}$, $\{S^4_k,S^5_k\}$ and $(\{C_i:w_s\in m_i\}\cup \{T_s\})$, $\{T_s,T^1_s\}$, $\{T^1_s,T^2_s\}$, $\{T^2_s, T^3_s\}$, $\{T^3_s,T^4_s\}$, $\{T^4_s,T^5_s\}$, respectively. Finally, we consider another clique $\{A_i,B_i,C_i:1\leq i\leq p\}$. Therefore, $\mathcal{T}$ is a tree with $8p+18q+1$ vertices. The construction of the clique tree is illustrated in Figure \ref{fig:Undirectedpath}.

	The graph $G$ is an undirected path graph having vertex set $\{A_i,B_i,C_i,D_i,\\E_i,F_i,G_i,H_i,J_i,K_i,O_i,P_i,Q_i,N_i,Y_i,Z_i:1\leq i\leq p\}\cup \{R_j,R^1_j,R^2_j,R^3_j,R^4_j,R^5_j,\\ S_k,S^1_k,S^2_k,S^3_k,s^4_k,S^5_k,T_s,T^1_s,T^2_s,T^3_s, T^4_s,T^5_s:1\leq j,k,s\leq q\}$ and two vertices of $G$ has an edge if they are in the same maximal clique. Clearly, $\mathcal{T}$ is a clique tree of $G$. Now, we prove the following claim.
	
	\begin{cl}\label{lem:UndirectedpathNPC}
		The graph $G$ has a liar's ve-dominating set of size $4p+10q$ if and only if there is a solution of \textsc{3-DMP}. 
	\end{cl}
	
	\begin{proof}
		Let $L$ be a liar's ve-dominating set of $G$ of size $4p+10q$. Consider the subgraph $G'$ induced by $\{A_i,B_i,C_i,D_i,E_i,F_i,G_i,H_i,J_i,K_i,P_i,O_i ,Q_i,N_i,Y_i,\\Z_i\}$ corresponding to the triple $m_i$. Note that, $L$ must contain at least four vertices of $G'$. Now if $L$ contains exactly four vertices of $G'$, then $L$ must contain the set $\{D_i,E_i,F_i,G_i\}$. Otherwise, $L$ contains five vertices of $G'$. If some vertices of $\{A_i,B_i,C_i,E_i,G_i\}$ are not in $L$, then we can obtain the set by replacing the missing vertex from this set. Thus, without loss of generality, we can assume that those five vertices are $\{A_i,B_i,C_i,E_i,G_i\}$. Let us assume that $L$ contains $\{A_i,B_i,C_i,E_i,G_i\}$ for $t$ number of $m_i$'s. Therefore, we have $|L|\geq 5t+4(p-t)$. Now, consider the subgraph induced by $\{R_j,R^1_j,R^2_j,R^3_j,R^4_j,R^5_j\}$. Corresponding to this subgraph $L$ must contain at least three vertices from the set $\{R^2_j,R^3_j,R^4_j,R^5_j\}$. Similarly corresponding to each subgraph induced by $ \{S_j,S^1_j,S^2_j,S^3_j,S^4_j,S^5_j\}$ and $ \{T_j,T^1_j,T^2_j,T^3_j,T^4_j,T^5_j\}$, $L$ must contain at least three vertices. Also, if $\{A_i,B_i,C_i\}\notin L$, then $L$ must contain at least $\max\{3(q-t),0\}$ number of $\{R_j,S_j,T_J\}$. Hence, we have $|L|\geq 5t+4(p-t)+3(q-t)+9q$. This implies that $t\geq q$. Thus, $L$ contains $q$ many $\{A_i,B_i,C_i\}$ vertices. Picking the corresponding $m_i$ forms a matching $M'$ of size $q$.   
		
		
		Conversely, let $M'$ be a solution of \textsc{3-DMP}. Consider the set $L=\{A_i,B_i,\\C_i,E_i,G_i:m_i\in M'\}\cup \{D_i,E_i,F_i,G_i:m_i\notin M'\}\cup \{R^2_j,R^3_j,R^4_j,S^2_k,S^3_k,S^4_k,T^2_s,\\T^3_s,T^4_s:1\leq j,k,s\leq q\}$. It is easy to see that $L$ is a liar's ve-dominating set of $G$ of cardinality $5q+4(p-q)+9q=4p+10q$.
	\end{proof}
	From the above lemma, we can conclude that \textsc{DecideLVEDP} is NP-complete for undirected path graphs.
\end{proof}

\section{Conclusion}
In this article, we have shown that the \textsc{MinLVEDP} can be solved in linear time for block graphs and proper interval graphs. We have proved that the \textsc{DecideLVEDP} is NP-complete for undirected path graphs. It is interesting to study the complexity problem for other subclasses of chordal graphs, bipartite graphs, and other graph classes. Also, the study of approximation algorithms for the \textsc{MinLVEDP} in different graph classes is needed. 
\bibliographystyle{plain}
\bibliography{VEDom_bib}

\begin{thebibliography}{10}

\bibitem{ahangar2021total}
H~Abdollahzadeh Ahangar, M~Chellali, SM~Sheikholeslami, M~Soroudi, and
  L~Volkmann.
\newblock Total vertex-edge domination in trees.
\newblock {\em Acta Mathematica Universitatis Comenianae}, 90(2):127--143,
  2021.

\bibitem{liar'svefirst}
Debojyoti Bhattacharya and Subhabrata Paul.
\newblock Algorithmic study on liar’s vertex-edge domination problem.
\newblock {\em Journal of Combinatorial Optimization}, 48(3):25, 2024.

\bibitem{Booth}
Kellogg~S. Booth and George~S. Lueker.
\newblock Testing for the consecutive ones property, interval graphs, and graph
  planarity using pq-tree algorithms.
\newblock {\em Journal of Computer and System Sciences}, 13(3):335--379, 1976.

\bibitem{boutrig2016vertex}
R.~Boutrig, M.~Chellali, T.~W Haynes, and S.T. Hedetniemi.
\newblock Vertex-edge domination in graphs.
\newblock {\em Aequationes mathematicae}, 90:355--366, 2016.

\bibitem{Totalve-domchellali}
Razika Boutrig and Mustapha Chellali.
\newblock Total vertex-edge domination.
\newblock {\em International Journal of Computer Mathematics},
  95(9):1820--1828, 2018.

\bibitem{chitra2012global}
S.~Chitra and R.~Sattanathan.
\newblock Global vertex-edge domination sets in graph.
\newblock In {\em International Mathematical Forum}, volume~7, pages 233--240,
  2012.

\bibitem{globalvedom}
S.~Chitra and R.~Sattanathan.
\newblock Global vertex-edge domination sets in total graph and product graph
  of path {$P_n$} cycle {$C_n$}.
\newblock In {\em Mathematical modelling and scientific computation}, volume
  283 of {\em Commun. Comput. Inf. Sci.}, pages 68--77. Springer, Heidelberg,
  2012.

\bibitem{garey1979}
M.R. Garey and D.S. Johnson.
\newblock {\em Computers and intractability:A Guide to the Theory of
  NP-Completeness}.
\newblock W.H Freeman, New York, 1979.

\bibitem{GAVRIL1975237}
Fǎnicǎ Gavril.
\newblock A recognition algorithm for the intersection graphs of directed paths
  in directed trees.
\newblock {\em Discrete Mathematics}, 13(3):237--249, 1975.

\bibitem{jena}
S.K. Jena and G.K. Das.
\newblock Vertex-edge domination in unit disk graphs.
\newblock {\em Discrete Applied Mathematics}, 319:351--361, 2022.

\bibitem{krishna}
B.~Krishnakumari, M.~Chellali, and Y.B. Venkatakrishnan.
\newblock Double vertex-edge domination.
\newblock {\em Discrete Mathematics, Algorithms and Applications},
  9(04):1750045, 2017.

\bibitem{lewis}
J.R. Lewis.
\newblock {\em Vertex-edge and edge-vertex parameters in graphs.}
\newblock PhD thesis, Clemson University, Clemson, SC, USA, 2007.

\bibitem{li2023polynomial}
Peng Li and Aifa Wang.
\newblock Polynomial time algorithm for k-vertex-edge dominating problem in
  interval graphs.
\newblock {\em Journal of Combinatorial Optimization}, 45(1):45, 2023.

\bibitem{naresh}
H.~Naresh~Kumar, D.~Pradhan, and Y.B. Venkatakrishnan.
\newblock Double vertex-edge domination in graphs: complexity and algorithms.
\newblock {\em Journal of Applied Mathematics and Computing}, 66(1):245--262,
  2021.

\bibitem{paul2}
S.~Paul, D.~Pradhan, and S.~Verma.
\newblock Vertex-edge domination in interval and bipartite permutation graphs.
\newblock {\em Discussiones Mathematicae: Graph Theory}, 43(4):947--963, 2021.

\bibitem{paul}
S.~Paul and K.~Ranjan.
\newblock Results on vertex-edge and independent vertex-edge domination.
\newblock {\em Journal of Combinatorial Optimization}, 44(1):303--330, 2022.

\bibitem{peters}
K.W.J. Peters.
\newblock {\em Theoritical and algorithmic results on domination and
  connectivity}.
\newblock PhD thesis, Clemson University, Clemson, SC, USA, 1986.

\bibitem{zylinski2019vertex}
P.~{\.Z}yli{\'n}ski.
\newblock Vertex-edge domination in graphs.
\newblock {\em Aequationes Mathematicae}, 93(4):735--742, 2019.

\end{thebibliography}

\end{document}